\documentclass[oneside, reqno, 11pt, a4paper]{amsart}

\usepackage{microtype}
\AtBeginDocument{
\DeclareSymbolFont{AMSb}{U}{msb}{m}{n}
\DeclareSymbolFontAlphabet{\mathbb}{AMSb}
}

\usepackage[dvipsnames]{xcolor}
\usepackage{graphicx}
\usepackage{tikz}
\usepackage{subfig}

\usepackage[a4paper, pdftex, left=2cm, top=2cm, right=2cm, bottom=2cm]{geometry}
\usepackage[foot]{amsaddr}

\usepackage{url}
\usepackage{hyperref}
\hypersetup{
	breaklinks=true,
	bookmarksopen=true,
	pdftitle={An intrinsic finite element framework for scalar- and vector-valued PDEs on general manifolds},
	pdfauthor={Tamara Tambyah, Alberto Martin, David Lee, Santiago Badia},
	pdfsubject={partial differential equations on manifolds, compatible finite elements},
	pdfkeywords={manifold, finite element exterior calculus, compatible finite elements, shallow water equations, cubed sphere, geometric error},
	colorlinks=true,
	linkcolor=black,
	citecolor=blue,
	filecolor=black,
	urlcolor=blue
}
\usepackage[capitalise]{cleveref}
\crefname{equation}{}{}
\crefrangeformat{equation}{(#3#1#4)--(#5#2#6)}

\usepackage{csquotes}
\usepackage[backend=biber, defernumbers=false, maxbibnames=99, style=numeric-comp, isbn=false, bibencoding=utf8, safeinputenc, url=false, doi=true, giveninits=true,sorting=none,natbib=true,uniquelist=false,maxcitenames=2]{biblatex}

\usepackage{float}
\usepackage{booktabs}
\usepackage{multirow}
\usepackage{algorithm}
\usepackage{algpseudocode}
\usepackage[inline]{enumitem}
\usepackage{acronym}

\acrodef{fe}[FE]{Finite Element}
\acrodefplural{fe}[FEs]{Finite Elements}
\acrodef{pde}[PDE]{Partial Differential Equation}
\acrodefplural{pde}[PDEs]{Partial Differential Equations}
\acrodef{ssprk3}[SSPRK3]{third order, strong stability preserving Runge--Kutta method}
\acrodef{api}[API]{Application Programming Interface}
\acrodef{jit}[JIT]{Just-In-Time}
\acrodef{supg}[SUPG]{streamline upwind Petrov--Galerkin} % define acros externally

\newtheorem{theorem}{Theorem}[section]

\newtheorem{proposition}[theorem]{Proposition}
\newtheorem{corollary}[theorem]{Corollary}
\newtheorem{definition}[theorem]{Definition}
\newtheorem{assumption}[theorem]{Assumption}
\newtheorem{remark}[theorem]{Remark}

\usepackage{amsmath, amsfonts, amssymb, amscd}
\usepackage{mathrsfs}
\usepackage{mathtools}
\usepackage[only, llbracket, rrbracket]{stmaryrd}
\allowdisplaybreaks
\usepackage[mathscr]{euscript}

\usepackage{lineno}

\usepackage{xspace}
\DeclareMathAlphabet{\mathpzc}{OT1}{pzc}{m}{it}

\DeclareMathOperator{\spn}{span}
\DeclareMathOperator{\codim}{codim}

\newcommand{\thesurface}{\mathcal{S}}
\newcommand{\surf}[1]{ \widetilde{{#1}} }
\newcommand{\surfvec}[1]{ \surf{\boldsymbol{{#1}}} }
\newcommand{\surfgrad}{ \nabla_{\thesurface} }
\newcommand{\surfdiv}{ \nabla_{\thesurface}\cdot }

\newcommand{\surfcurl}{ \nabla_{\thesurface} \times }

\newcommand{\V}{\mathbb{V}}

\newcommand{\gravity}{ \mathpzc{g} } 
\newcommand{\coriolis}{ \mathpzc{f} } 
\newcommand{\topography}{ \mathpzc{b} } 

\newcommand{\map}{\sigma}
\newcommand{\pull}{\map^*}

\newcommand{\thechart}{\mathcal{V}}
\newcommand{\chart}[1]{ {#1} }

\newcommand{\RT}{Raviart-Thomas\xspace}
\newcommand{\N}{N\'ed\'elec\xspace}

\newcommand{\twoD}{two-dimensional\xspace}
\newcommand{\threeD}{three-dimensional\xspace}

\newcommand{\radius}{\mathscr{R}}
\newcommand{\thickness}{\mathscr{T}}

\newcommand{\GridapGeoscience}{GridapGeosciences.jl\xspace}

\newcommand{\gridap}{Gridap\xspace}
\newcommand{\GridapDistributed}{GridapDistributed.jl\xspace}
\newcommand{\GridapSolvers}{GridapSolvers.jl\xspace}
\newcommand{\GridapPETSc}{GridapPETSc.jl\xspace}
\newcommand{\GridapPest}{GridapP4est.jl\xspace}

\newcommand{\Riemannian}{Riemannian\xspace}

\newcommand{\discretiz}{discretis\xspace}
\newcommand{\discretization}{{\discretiz}ation\xspace}

\newcommand{\discretizations}{{\discretiz}ations\xspace}

\newcommand{\RK}{Runge--Kutta\xspace}

\newcommand{\TS}{ T_{\surf{x}}\thesurface }

\newcommand{\contra}[1]{\vec{#1}}            % contravariant component field
\newcommand{\covar}[1]{\underline{#1}}       % covariant component field
\newcommand{\gradg}{\nabla_{g}}              % intrinsic gradient
\newcommand{\divg}{\nabla_{g}\cdot}          % intrinsic divergence
\newcommand{\curlg}{\nabla_{g}\times}        % intrinsic curl (n=3)
\newcommand{\lapg}{\Delta_{g}}               % Laplace--Beltrami
\newcommand{\gradgperp}{\nabla_{g}^{\perp}}  % intrinsic skew gradient (n=2)
\newcommand{\divgperp}{\nabla_{g}^{\perp}\cdot} % intrinsic skew divergence (n=2)
\newcommand{\transm}{\mathsf{T}}             % transmission map between charts
\newcommand{\flux}[1]{\boldsymbol{#1}}       % flux proxy (density representation): sqrt(g) x contravariant

\let\oldalign\align
\let\oldendalign\endalign
\renewenvironment{align}
{\linenomathNonumbers\oldalign}
{\oldendalign\endlinenomath}

\title[Intrinsic finite elements on manifolds]{
An intrinsic finite element framework for scalar- and vector-valued partial differential equations on general manifolds}
\date{\today}
\keywords{manifold, finite element exterior calculus, compatible finite elements, shallow water equations, cubed sphere, geometric error}

\author[1]{Tamara A. Tambyah$^{1}$}
\address{$^1$School of Mathematics\\Monash University\\Clayton\\Victoria 3800\\Australia }
\email{tamara.tambyah@monash.edu}

\author[2]{Alberto F. Mart\'in$^{2}$}
\email{alberto.f.martin@anu.edu.au}
\address{$^{2}$School of Computing\\Australian National University\\Canberra\\ACT 2600\\Australia }

\author[3]{David Lee$^{3}$}
\email{david.lee@bom.gov.au}
\address{ $^{3}$Bureau of Meteorology\\Melbourne\\Australia}

\author[4]{Santiago Badia$^{1\ast}$}
\email{santiago.badia@monash.edu}

\thanks{$^\ast$Corresponding author}

\AtEveryCitekey{\clearfield{eprint}} % remove eprints in bibliography

\usepackage[normalem]{ulem}

\begin{document}

\begin{abstract}
We present an intrinsic finite element framework for the numerical approximation of scalar- and vector-valued partial differential equations on general manifolds described by atlases of charts.
Weak formulations and their \discretization are expressed exclusively in the flat parametric space of the manifold, where the exact geometry enters only through the metric tensor that is evaluated at quadrature points. 
Conforming finite element spaces of arbitrary order for the whole de Rham complex are obtained on multi-chart atlases through purely topological inter-chart gluing.
We develop a nodal change of basis construction for vector-valued spaces on multi-chart atlases, where nodal degrees of freedom couple charts with incompatible coordinate systems through transmission maps. 
Inter-chart continuity holds at the nodes with jumps of the order of the approximation error, and the resulting discrete field is exactly tangent to the manifold. Optimal convergence rates for this method are observed numerically for the surface Stokes problem; no a priori error analysis is provided. 
Numerical experiments on the cubed sphere manifold in two and three dimensions for the Hodge Laplacian problems
demonstrate the intrinsic framework is free of the geometric consistency errors for a sufficient quadrature degree, and that intrinsic assembly is substantially cheaper than its extrinsic counterpart. 
Further considering the rotating shallow water equations that include  orography as a perturbation of the geometry demonstrates that mass is conserved exactly, and energy is conserved up to the time \discretization error.

\end{abstract}

\maketitle

\section{Introduction}
\label{sec: introduction}

Many applications involve \acfp{pde} posed on manifolds. 
A prominent example is numerical weather prediction, where dynamical cores solve \acp{pde} for geophysical fluid dynamics on the sphere or a thin atmospheric shell \cite{Staniforth2012,Melvin2024_orography}.
In this context, compatible \acf{fe} methods are an established \discretization strategy \cite{Cotter2012,Mcrae2013,Natale2016,Cotter2023,BauerCotter2018}, 
where  \ac{fe} spaces form a discrete de Rham complex \cite{Arnold2006,Arnold2018} that inherits vector calculus identities at the discrete level. 
This yields local mass conservation, steady geostrophic modes, and suppresses spurious modes.
A \discretization on a manifold must therefore respect two structures:
(1) the geometry, and (2) the differential structure of the equations, where 
unknowns satisfy mimetic identities, and 
vector-valued fields are tangent to the manifold.
This work presents a generic \ac{fe} framework that exactly preserves the differential structure at any \discretization order, and the geometry up to numerical quadrature errors.

The prevailing \ac{fe} treatment of \acp{pde} on manifolds is \emph{extrinsic}. 
That is, a smooth manifold
is approximated via a piecewise-polynomial surface triangulation, and the weak formulation is posed in ambient coordinates \cite{Dziuk2013,Bonito2019}.
Approximating the manifold introduces \emph{geometric consistency error} that adds to the approximation error of the \ac{fe} solution.  
Such a variational crime \cite{HolstStern2012} dominates the convergence of the method \cite{Dziuk2013,Bonito2019,Holst2023}, 
unless the order of the geometry representation increases with that of the \ac{fe} spaces \cite{Stanford2019,Fries2017}.  
This extrinsic viewpoint underpins many \ac{fe} software stacks that support immersed manifolds, such as  FEniCS \cite{Rognes2013} and Firedrake \cite{Gibson2019_firedrake}, as well as  the compatible \ac{fe} dynamical core of the UK Met Office next-generation atmospheric model, LFRic, where the cubed sphere manifold represents the Earth \cite{Adams2019,Kent2023}. 
In such a weather model,
atmospheric dynamics are resolved using
mixed \ac{fe} spaces that map via Piola transformations onto a polynomial approximation of the cubed sphere  \cite{Melvin2019,Melvin2024_orography}.
The resulting geometry is $C^{0,1}$,
where tangent planes are not well defined across the element boundaries.
Piola transformations only guarantee element-wise tangency to the approximate surface, 
meaning the resulting \discretization is non-conforming in general, and geometrical consistency error with respect to the true sphere  persists.
These variational crimes are eliminated in parametric methods that utilise exact analytical maps to describe the manifold \cite{Guba2014,Taylor2010}. 
However, analytical descriptions cannot be readily accommodated in general \ac{fe} software stacks.

Extrinsic surface \ac{fe} methods
are well-developed for scalar unknowns, with a complete a priori error theory \cite{Dziuk2013,Bonito2019}. 
These methods employ a piecewise-polynomial surface that is  a $C^{0,1}$ geometry, for which  vector-valued spaces are well defined in element interiors. 
However, pointwise tangency of continuous vector-valued fields cannot be imposed on both sides of inter-element interfaces  since the tangent spaces of adjacent edge elements differ.  
Consequently, tangency of the discrete vector field and $H^1$-conformity of the \ac{fe} space cannot be strongly enforced simultaneously.
In many cases, tangency is imposed weakly, by penalising the normal component of a fully ambient velocity \cite{Hansbo2020}, or via Lagrange multipliers and trace \ac{fe} techniques \cite{Gross2018}. 
Alternatively, exact tangency can be retained by abandoning $H^1$-conformity altogether. 
In such cases, continuity is weakly restored through discontinuous Galerkin couplings of Piola-mapped elements \cite{Lederer2020}.
Extrinsic surface \ac{fe} methods ultimately remain problematic for vector-valued unknowns    \cite{DemlowNeilan2024,Hansbo2020,DemlowNeilan2026}, 
and demand a dedicated analysis of the tangential accuracy of the discrete solution \cite{DemlowNeilan2026,Hansbo2020,Hardering2023}.

An \emph{intrinsic} viewpoint alleviates the incompatibility between tangency and $H^1$-conformity of vector-valued fields, by posing the weak formulation in the flat parametric space of the manifold, and encoding the geometry through the \Riemannian metric tensor.
For scalar unknowns,
\citet{Bachini2021} proposed an intrinsic \ac{fe} method for  advection--diffusion--reaction equations with lowest-order elements, and surfaces described by a single chart;
the extension to multiple charts is there identified as a limitation.
The same single-chart setting underlies the intrinsic derivation of the shallow water equations developed by \citet{Bachini2020}.
Hence, to the best of our knowledge, the literature lacks an intrinsic \emph{compatible} \ac{fe} framework that supports the full de Rham complex and vector-valued unknowns on manifolds described by multi-chart atlases.

To fill this gap, we develop a generic intrinsic \ac{fe} framework in which both the weak formulation and its \ac{fe} \discretization\ are expressed exclusively in the parametric space of an atlas, where the exact geometry enters only through the metric tensor at numerical integration points.
This yields a geometry that is exact and smooth, meaning tangency and $H^1$-conformity coexist within each chart of the manifold atlas. 
The incompatibility arises only at chart interfaces with mismatched coordinate frames, for which we enforce inter-chart continuity at the nodes via  a novel construction in contribution (ii) below.
We deliberately adopt vector calculus--based approach as opposed to an  exterior calculus one \cite{Cotter2014} so that our presentation directly maps to \ac{fe} codes, most of which rely on data structures originating from vector calculus notation.
The specific contributions of this study are:
\begin{enumerate}[label=(\roman*)]
	\item We derive fully intrinsic weak formulations for scalar- and vector-valued \acp{pde} on multi-chart atlases
	using a vector calculus exposition amenable to practitioners.
	Such an approach is also used to derive the mimetic identities and the skew operators central to atmospheric models. 
	\item We construct parametric \ac{fe} spaces of arbitrary order for the whole discrete de Rham complex where tangency of the resulting fields holds by construction. 
	For grad-conforming vector-valued fields, whose degrees of freedom couple charts with incompatible coordinate systems, we introduce a novel nodal change of basis at chart interfaces, where inter-chart continuity holds at the nodes and the interface jumps are of the order of the approximation error.
	The resulting fields are exactly tangent;  rigorous convergence analysis is left as an open question.
 
	\item 
	We demonstrate intrinsic formulations are free of geometric consistency errors,
	since the only geometry-related error is the numerical approximation of non-polynomial metric terms, which reduces as the quadrature degree increases.
	
	\item We quantify the computational benefits of an intrinsic pipeline, for which the metric is supplied analytically at quadrature points and all  computations are performed in the parametric space, in comparison to  an   extrinsic method that utilises the ambient space of the manifold.  
	\item We implement the intrinsic framework in \GridapGeoscience \cite{GridapGeo,Tambyah2026_GridapGeo}, which utilises the open-source  \gridap \ac{fe} ecosystem  in Julia \cite{Badia2020,Verdugo2022} and supports distributed-memory parallelism. 
\end{enumerate}

Several model problems are used to showcase the intrinsic framework on the cubed sphere manifold \cite{Ronchi1996}, for which the atlas of six charts is used widely in atmospheric modelling \cite{Staniforth2012,Putman2007}.
First, we solve the Hodge Laplacian problems associated to the discrete de Rham complex \cite{Arnold2018},
and demonstrate the non-polynomial cubed sphere metric is captured to machine precision for a sufficient quadrature degree, with the expected $hp$-convergence.
We then discretise the rotating shallow water equations, which is the standard \twoD testbed for atmospheric dynamical cores \cite{Williamson1992,Mcrae2013,Cotter2014}. 
We obtain a compatible formulation with consistent potential vorticity upwinding \cite{Lee2024},
for which mass conservation is exact and energy conservation holds up to the time \discretization error. 
The Williamson test suite \cite{Williamson1992} confirms the expected convergence rates and conservation properties.
The intrinsic framework is agnostic to the choice of metric, of which many variants exist for the cubed sphere alone \cite{Rancic1996,Giraldo2003,Nair2005,Mcgregor2005}.
We exploit this flexibility to solve the shallow water equations on a \threeD atmospheric shell with a single radial element, 
where orography is included as a perturbation of the geometry itself, 
as a stepping stone towards terrain-following atmospheric models \cite{Wood2013,Melvin2024_orography}. 
Finally, 
{we observe optimal convergence rates in a discontinuous Galerkin-like norm \cite{Cockburn2009} when applying the developed nodal construction of grad-conforming vector fields to the surface Stokes problem on the cubed sphere manifold.}

The remainder of this article is structured as follows.
\Cref{sec: math form} introduces the parametric representation of manifolds, and the intrinsic form of the differential operators and mimetic identities.
\Cref{sec: fem} constructs the parametric \ac{fe} spaces and discrete de Rham complexes, including the grad-conforming vector-valued elements on multi-chart atlases.
\Cref{sec: models} derives the intrinsic formulations of the Hodge Laplacian problems and the rotating shallow water equations, in addition to their discrete conservation properties.
\Cref{sec: results} presents the numerical experiments, and conclusions are drawn in \cref{sec: conclusion}.

\section{Mathematical notation}
\label{sec: math form}

This section establishes the notation used throughout the article.
We deliberately restrict the presentation to the minimal set of geometric objects that the intrinsic framework requires: an atlas of parametric domains and geometrical maps, the representation of tangent fields by their component fields, and the metric tensor.
Where appropriate, we relate concepts to the cubed sphere manifold.
The reader is referred to \citet{Frankel2004} and \citet{Lee2003_manifold} for a comprehensive introduction to smooth manifolds.

%%%%%%%%%%%%%%%%%%%%%%%%%%%%%%%%%%%%%%%%%%%%%%%%%%%%%%%%%%%%%%%%%%%%%%%%%%%%%%
\subsection{Manifolds, charts and fields}
\label{sec: manifolds}

Let $\thesurface \subset \mathbb{R}^m$ be a smooth, orientable, compact $n$-dimensional manifold, where $n \leq m$ is the topological dimension, $m$ is the ambient dimension and $\codim(\thesurface) = m - n$ is the co-dimension.
We describe $\thesurface$ by an \emph{atlas} $\left\{ (\thechart_k, \map_k) \right\}_{k=1}^{\kappa}$ of \emph{charts} $(\thechart_k,\map_k)$. 
The \emph{parametric domains} $\thechart_k \subset \mathbb{R}^n$ are open and bounded, the \emph{geometrical maps} $\map_k : \thechart_k \rightarrow \mathbb{R}^m$ are smooth injective immersions, and the union of the images $\map_k(\thechart_k)$ covers $\thesurface$.
{We further assume a smooth embedding such that $\map_k$ is a homeomorphism onto its image.}

Where two images overlap, quantities on $\thesurface$ can be described by two different charts $k$ and $l$. 
In this instance, the smooth \textit{transition map} $\tau_{kl} = \map_l^{-1} \circ \map_k$  converts points in chart $k$ and $l$ of the overlap into the same point on $\thesurface$. 
That is, $\map_k(\chart{x}) = \map_l(\tau_{kl}(\chart{x}))$
such that chart-wise quantities are consistent on overlaps.
{The atlas is assumed to be \textit{oriented} such that  $\det \mathrm{d}\tau_{kl} > 0$.}
Ambient quantities carry a tilde and their parametric counterparts are plain such that $\chart{x} = (x^1, \dots, x^n)$ denotes a point in $\thechart_k$, and $\surf{x} = \map_k(\chart{x})$ is the corresponding point on $\thesurface$.
Since all definitions below are stated chart-wise, we work with a single chart $(\thechart, \map)$ and drop the index $k$.

Given a scalar field $\surf{f} : \thesurface \rightarrow \mathbb{R}$, its \emph{pullback} by the geometrical map is
\begin{align}
	f = \pull \surf{f} = \surf{f} \circ \map ,
	\label{eq: pullback}
\end{align}
where $f : \thechart \rightarrow \mathbb{R}$. 
We systematically denote a field on $\thesurface$ and its pullback by the same letter, with and without a tilde, respectively.  
That is,  a tilde identifies fields on $\thesurface$, plain symbols denote parametric quantities and \cref{eq: pullback} is the relation between the two. 
The same convention applies to vector quantities.
We write $\surfvec{v}$ for a tangent vector field on $\thesurface$, and $\contra{v}$ and $\covar{v}$ for its contravariant and covariant component fields on $\thechart$, whose precise definitions follow below.
Boldface with a tilde is reserved for $\mathbb{R}^m$-valued fields on $\thesurface$, arrows and underlines identify $\mathbb{R}^n$-valued parametric component fields, and boldface without a tilde is reserved for the flux proxies introduced in \cref{sec: pullbacks} below.
Ambient vector fields are confined to this section and to the statement of model problems on $\thesurface$. 
All developments that follow involve parametric objects alone.

Since $\map$ is an immersion, its differential $\mathrm{d}\map_{\chart{x}} \in \mathbb{R}^{m \times n}$, evaluated at every $\chart{x} \in \thechart$,  has full rank $n$.
The \emph{tangent space} of $\thesurface$ at $\surf{x} = \map(\chart{x})$ is the $n$-dimensional subspace of $\mathbb{R}^m$
\begin{align}
	\TS = \mathrm{range} \left( \mathrm{d}\map_{\chart{x}} \right)
	= \spn \left\{ \surf{\boldsymbol{\partial}}_1(\chart{x}), \dots, \surf{\boldsymbol{\partial}}_n(\chart{x}) \right\},
	\qquad
	\text{where},
	\qquad
	\surf{\boldsymbol{\partial}}_i = \frac{\partial \map}{\partial x^i} ,
	\label{eq: tangent space}
\end{align}
are the \emph{covariant basis vectors} $\surf{\boldsymbol{\partial}}_i : \thechart \rightarrow \mathbb{R}^m$ that form the columns of $\mathrm{d}\map$.
A \emph{tangent vector field} is a map $\surfvec{v} : \thesurface \rightarrow \mathbb{R}^m$ such that $\surfvec{v}(\surf{x}) \in \TS$ at every point.
Since $\mathrm{d}\map$ has full rank, every tangent vector field is uniquely represented by its \emph{contravariant component field} $\contra{v} : \thechart \rightarrow \mathbb{R}^n$, defined by
\begin{align}
	\surfvec{v} \circ \map = \mathrm{d}\map \, \contra{v},	
	\label{eq: contravariant expansion}
\end{align}
In components,  $\surfvec{v} = v^i \surf{\boldsymbol{\partial}}_i$ with summation over repeated indices.

The \emph{metric tensor} $g : \thechart \rightarrow \mathbb{R}^{n \times n}$ is defined by
\begin{align}
	g = \mathrm{d}\map^T \mathrm{d}\map ,
	\qquad \text{with components}, \quad
	g_{ij} = \surf{\boldsymbol{\partial}}_i \cdot \surf{\boldsymbol{\partial}}_j ,
	\label{eq: metric}
\end{align}  
and is a symmetric positive definite matrix at every $\chart{x} \in \thechart$, since $\map$ is an immersion.
We write $g^{-1}$ with components $g^{ij}$ for the metric inverse and $\sqrt{g} = (\det g)^{1/2}$ for the square root of the metric determinant.
The metric encodes the pointwise inner product of tangent fields,
\begin{align}
	\left( \surfvec{v} \cdot \surfvec{w} \right) \circ \map = \contra{v} \cdot g \, \contra{w} ,
	\label{eq: inner product}
\end{align}
and defines the \emph{covariant component field} of $\surfvec{v}$ as $\covar{v} = g \, \contra{v}$, which in components is $v_i = g_{ij}v^j$ such that 
$g$ lowers the contravariant index. 
Conversely, $\contra{v} = g^{-1} \covar{v}$ 
such that $g^{-1}$ raises the covariant index as $v^i = g^{ij}v_j$. 
We stress that a tangent field $\surfvec{v}$ is an invariant object. The adjectives covariant and contravariant refer to how its components transform under a change of coordinates, not to the tangent field itself.
Finally, the metric relates integration on $\thesurface$ and on the parametric domain:
\begin{align}
	\int_{\thesurface} \surf{f} \, \mathrm{d}\thesurface
	= \int_{\thechart} f \sqrt{g} \, \mathrm{d}\chart{x} ,
	\label{eq: integration}
\end{align}
for scalar fields supported on $\map(\thechart)$, 
{where $\mathrm{d}\thesurface$ is the surface measure on the manifold $\thesurface$ and $\mathrm{d}\chart{x}$ the Lebesgue measure in the parametric domain $\thechart$.}

To integrate over the whole manifold in the multi-chart case, we introduce a \emph{partition subordinate to the atlas}. That is, closed subdomains $P_k \subset \thechart_k$ whose images $\map_k(P_k)$ cover $\thesurface$ and pairwise intersect on sets of measure zero, so that integrals over $\thesurface$ are the sum of the contributions of the partition domains. 

\begin{remark}
	\label{rem: intrinsic}
	The framework is \emph{intrinsic} since 
	the subsequent development of differential operators, weak formulations and their \ac{fe} \discretization involves only the metric $g$, together with $g^{-1}$ and $\sqrt{g}$.
	Once the metric is provided, weak formulations and spaces of the de Rham complex can be constructed, while grad-conforming spaces, described in \cref{sec: grad conforming} below, require the transition maps of the atlas. 
	For many manifolds of practical interest, such as the cubed sphere \cite{Ronchi1996}, the metric is available in closed form and can be evaluated analytically at quadrature points, without ever forming $\mathrm{d}\map$ or computing in the ambient dimension $m$.
	In fact, the geometrical map $\map$ and its differential $\mathrm{d}\map$ are only appear in \cref{eq: tangent space,eq: metric,eq: contravariant expansion} above, and are used 	
	to pull back prescribed ambient data in \cref{sec: models} below.
\end{remark}

\cref{tab: notation} summarises the notation introduced so far, together with the intrinsic and surface differential operators defined in \cref{sec: intrinsic operators} below.

\begin{table}[h!]
	\centering
	\footnotesize
	\setlength{\tabcolsep}{3pt}
	\resizebox{\ifdim\width>\textwidth\textwidth\else\width\fi}{!}{%
	\begin{tabular}{ll@{\quad}ll}
		\toprule
		Symbol & Meaning & Symbol & Meaning \\
		\midrule
		$\thesurface$, $n$, $m$ & manifold, topological/ambient dimensions &
		$g$, $g^{-1}$, $\sqrt{g}$ & metric \cref{eq: metric}, inverse, $(\det g)^{1/2}$ \\
		$\thechart$, $\map$ & chart: parametric domain, geometrical map &
		$\nabla$, $\nabla \cdot$, $\nabla \times$ & Euclidean grad, div, curl ($n=3$) in $\thechart$ \\
		$P_k \subset \thechart_k$ & partition domains subordinate to the atlas &
		$\gradg$, $\divg$, $\lapg$ & intrinsic grad, div, Laplace--Beltrami \cref{eq: intrinsic ops} \\
		$\chart{x}$, $\surf{x} = \map(\chart{x})$ & parametric point, ambient point &
		$\curlg$ & intrinsic curl for $n = 3$ \cref{eq: intrinsic curl} \\
		$\surf{f}$, $f = \surf{f} \circ \map$ & scalar field on $\thesurface$, pullback on $\thechart$ &
		$\gradgperp$, $\divgperp$ & intrinsic skew grad, skew div for $n = 2$ \cref{eq: skew ops} \\
		$\surf{\boldsymbol{\partial}}_i$, $\TS$ & covariant basis \cref{eq: tangent space}, tangent space &
		$\pull_0, \pull_1, \dots$ & pullbacks \cref{eq: pullbacks 3D,eq: pullbacks 2D} \\
		$\surfvec{v} = v^i \surf{\boldsymbol{\partial}}_i$ & tangent vector field on $\thesurface$ \cref{eq: contravariant expansion} &
		$\covar{n}$ & Euclidean unit normal (covariant) to $\partial \thechart$  \\
		$\contra{v} = (v^i)$, $\covar{v} = g \contra{v}$ & contra-, covariant component fields &
		$\surfgrad$, $\surfdiv$, $\surfcurl$ & surface grad, div, curl ($n = 3$) (Def.~\ref{rem: ambient operators}) \\
		$\flux{v} = \sqrt{g} \, \contra{v}$ & flux proxy  &
		$\surfgrad^{\perp}$, $\surfgrad^{\perp} \cdot$ & surface skew grad, skew div for $n = 2$ (Def.~\ref{rem: ambient operators}) \\
		\bottomrule
	\end{tabular}}
	\vspace{0.5em}
	\caption{Summary of notation.}
	\label{tab: notation}
\end{table}

%%%%%%%%%%%%%%%%%%%%%%%%%%%%%%%%%%%%%%%%%%%%%%%%%%%%%%%%%%%%%%%%%%%%%%%%%%%%%%

\subsection{Intrinsic differential operators}
\label{sec: intrinsic operators}

Differential operators in the parametric domain are denoted  $\nabla = (\partial  / \partial x^1, \dots, \partial  / \partial x^n)$ for the Euclidean gradient, $\nabla \cdot$ for the Euclidean divergence, and  $\nabla \times$ for the Euclidean curl when $n=3$.
Intrinsic differential operators,
denoted with a subscript $g$,  act on and return \emph{parametric} fields. 
The corresponding tangent fields on $\thesurface$ are recovered through \cref{eq: contravariant expansion}.
For a scalar field $f$ and a tangent field with contravariant components $\contra{v}$, the intrinsic gradient, divergence and Laplace--Beltrami operators are, respectively, \cite{Frankel2004,Bachini2021}
\begin{align}
	\gradg f = g^{-1} \nabla f ,
	\qquad
	\divg \contra{v} = \frac{1}{\sqrt{g}} \nabla \cdot \left( \sqrt{g} \, \contra{v} \right) ,
	\qquad
	\lapg f = \divg \left( g^{-1} \nabla f \right) ,
	\label{eq: intrinsic ops}
\end{align}
where $\divg$ and $\lapg$ return scalar fields.
Note that $\nabla f=(\partial f / \partial x^1, \dots, \partial f / \partial x^n)$ collects the derivatives $\partial_i f$, which transform as covariant components. Raising the index with $g^{-1}$ means the intrinsic gradient $\gradg f$ returns contravariant components.

When $n = 3$, the intrinsic curl acts on the covariant components of a tangent field and returns contravariant components: 
\begin{align}
	\curlg \covar{v} = \frac{1}{\sqrt{g}} \nabla \times \covar{v} .
	\label{eq: intrinsic curl}
\end{align}

When $n = 2$, we introduce the rotation matrix, $R$, and the rotated Euclidean gradient
\begin{align*}
	\nabla^{\perp} = R \nabla = \left( - \frac{\partial}{\partial x^2}, \frac{\partial}{\partial x^1} \right),
	\qquad
	\text{where}
	\qquad
		R = \begin{pmatrix} 0 & -1 \\ 1 & 0 \end{pmatrix} .
\end{align*}
Then, the
intrinsic skew gradient and skew divergence operators are 
\begin{align}
	\gradgperp f = \frac{1}{\sqrt{g}} \nabla^{\perp} f ,
	\qquad
	\divgperp \contra{v} = - \frac{1}{\sqrt{g}} \nabla \cdot \left( R \, \covar{v} \right)
	= \frac{1}{\sqrt{g}} \left( \frac{\partial v_2}{\partial x^1} - \frac{\partial v_1}{\partial x^2} \right) ,
	\label{eq: skew ops}
\end{align}
where $\gradgperp f$ returns contravariant components and $\covar{v} = g \contra{v}$.
The operators in \cref{eq: skew ops} are the two-dimensional counterparts of the intrinsic curl \cref{eq: intrinsic curl}, where  $\gradgperp$ maps scalars to tangent fields, and $\divgperp$ maps tangent fields to scalars.
They are the fundamental building blocks of the vorticity dynamics in the shallow water equations of \cref{sec: models} below.

The intrinsic operators satisfy the following annihilation  identities, which are the backbone of the compatible \discretizations\ developed in this work.

\begin{proposition}
	\label{prop: annihilation}
	Let $f$ be a smooth scalar field and $\covar{v}$ a smooth covariant component field on $\thechart$. Then
	\begin{align}
		\curlg \left( \covar{\gradg f} \right) = \boldsymbol{0} ,
		\qquad
		\divg \left( \curlg \covar{v} \right) = 0 ,
		\qquad (n = 3),
		\label{eq: annihilation 3D}
		\\
		\divgperp \left( \gradg f \right) = 0 ,
		\qquad \qquad \qquad \quad ~
		(n = 2),
		\label{eq: annihilation 2D}
	\end{align}
	where $\covar{\gradg f} = g \gradg f = \nabla f$ denotes the covariant components of the intrinsic gradient.
	\begin{proof}
		In each identity, the metric contributions cancel exactly, and the result reduces to the Euclidean annihilation properties in $\thechart$ due to the symmetry of second derivatives.
		For the first identity, $\curlg (\nabla f) = \frac{1}{\sqrt{g}} \nabla \times (\nabla f) = \boldsymbol{0}$.
		For the second, $\divg ( \curlg \covar{v} ) = \frac{1}{\sqrt{g}} \nabla \cdot \left( \sqrt{g} \frac{1}{\sqrt{g}} \nabla \times \covar{v} \right) = \frac{1}{\sqrt{g}} \nabla \cdot ( \nabla \times \covar{v} ) = 0$.
		For the third, $\divgperp ( \gradg f ) = - \frac{1}{\sqrt{g}} \nabla \cdot ( R g \, g^{-1} \nabla f ) = - \frac{1}{\sqrt{g}} \nabla \cdot ( R \nabla f ) = 0$.
	\end{proof}
\end{proposition}

\begin{remark}
	\label{rem: exact annihilation}
	The identities in \cref{prop: annihilation} hold pointwise and for \emph{any} metric since the geometric-dependent terms cancel identically.
	As a consequence, the identities in \cref{prop: annihilation} are inherited exactly by the discrete spaces formed in \cref{sec: discrete complex} below. 
\end{remark}

\begin{definition}
	\label{rem: ambient operators}
	Surface differential operators, denoted with a subscript $\thesurface$, act on ambient fields on $\thesurface$.
	The following definitions coincide with the classical tangential operators, defined extrinsically through projections of ambient operators onto the tangent space or cross products with the unit normal \cite{Dziuk2013}, and are therefore independent of the chart.
	The \emph{surface gradient} $\surfgrad \surf{f}$ is the tangent field whose contravariant component field is $\gradg f$. The \emph{surface divergence} $\surfdiv \surfvec{v}$ is the scalar field whose pullback is $\divg \contra{v}$. 
	 For $n=3$, the \emph{surface curl}  $\surfcurl \surfvec{v}$ is induced by \cref{eq: intrinsic curl} where the intrinsic curl acts through the covariant component field.
	Likewise,  $\surfgrad^{\perp} \surf{f}$ and $\surfgrad^{\perp} \cdot \, \surfvec{v}$ arise from  \cref{eq: skew ops} for $n=2$.
\end{definition}
In this work, the surface operators in \cref{rem: ambient operators} are never used in computations; they appear only in the statement of problems and identities on $\thesurface$. 
No ambient machinery, in particular the normal vector to $\thesurface$, is required in the developments that follow.

%%%%%%%%%%%%%%%%%%%%%%%%%%%%%%%%%%%%%%%%%%%%%%%%%%%%%%%%%%%%%%%%%%%%%%%%%%%%%%
\subsection{Pullbacks, function spaces and weak identities}
\label{sec: pullbacks}

For each conformity class of fields on $\thesurface$,
we now present  the parametric representation of the differential operators defined in \cref{sec: intrinsic operators} above, and their integration by parts identities.
Throughout this section we assume that $\map$ extends to a smooth immersion on the closure of $\thechart$, so that the eigenvalues of $g$ are uniformly bounded from above and away from zero.

For $n = 3$, we define the \emph{pullbacks} of scalar and tangent vector fields as 
\begin{align}
	\pull_0 \surf{f} = f ,
	\qquad
	\pull_1 \surfvec{v} = \covar{v} ,
	\qquad
	\pull_2 \surfvec{v} = \flux{v} ,
	\qquad
	\pull_3 \surf{\rho} = \sqrt{g} \, \rho ,
	\label{eq: pullbacks 3D}
\end{align}
where the subscript refers to the degree of the differential form that each object proxies \cite{Arnold2018,HolstStern2012},
and $f$ and $\rho$ denote the scalar pullbacks \cref{eq: pullback} of $\surf{f}$ and $\surf{\rho}$.
In particular, $\pull_0$ coincides with the pullback $\pull$ of \cref{eq: pullback}.
The maps $\pull_1$ and $\pull_2$ are the covariant and contravariant Piola vector transformations \cite[Chapter 9]{ErnGuermond_1}.
In \cref{eq: pullbacks 3D}, the boldface field $\flux{v} = \sqrt{g} \, \contra{v}$ is the \emph{flux proxy} of $\surfvec{v}$, for which the contravariant components are scaled by $\sqrt{g}$.
In the language of \ac{fe} exterior calculus, the flux proxy, and not $\contra{v}$, is the object approximated by div-conforming \ac{fe} spaces \cite{Arnold2018}.
Boldface without a tilde is reserved for flux proxies hereinafter.
The pullbacks in \cref{eq: pullbacks 3D} define chart-wise
Sobolev spaces on $\thesurface$:
\begin{align}
	\surf{f} \in H^1(\thesurface)
	& \; \Longleftrightarrow \;
	\pull_0 \surf{f} \in H^1(\thechart) ,
	&
	\surfvec{v} \in H(\mathrm{curl}, \thesurface)
	& \; \Longleftrightarrow \;
	\pull_1 \surfvec{v} \in H(\mathrm{curl}, \thechart) ,
	\nonumber
	\\
	\surfvec{v} \in H(\mathrm{div}, \thesurface)
	& \; \Longleftrightarrow \;
	\pull_2 \surfvec{v} \in H(\mathrm{div}, \thechart) ,
	&
	\surf{\rho} \in L^2(\thesurface)
	& \; \Longleftrightarrow \;
	\pull_3 \surf{\rho} \in L^2(\thechart) ,
	\label{eq: conformity}
\end{align}
for every chart of the atlas.
By the uniform bounds on the eigenvalues of $g$, these definitions are equivalent to the usual ambient ones \cite{ErnGuermond_1,Dziuk2013,HolstStern2012}, and the parametric $L^2$ norms are equivalent to their $g$-weighted counterparts.

The pullbacks \cref{eq: pullbacks 3D,} intertwine the surface differential operators with the Euclidean operators on $\thechart$, as summarised in the following commuting diagram.
\begin{proposition}
	\label{prop: commuting}
	For smooth fields and $n = 3$,
	\begin{align}
		\pull_1 \left( \surfgrad \surf{f} \right) = \nabla \left( \pull_0 \surf{f} \right) ,
		\qquad
		\pull_2 \left( \surfcurl \surfvec{v} \right) = \nabla \times \left( \pull_1 \surfvec{v} \right) ,
		\qquad
		\pull_3 \left( \surfdiv \surfvec{v} \right) = \nabla \cdot \left( \pull_2 \surfvec{v} \right) ,
		\label{eq: commuting 3D}
	\end{align}
such that the following diagram commutes:
	\begin{equation}
		\begin{CD}
			H^1(\thesurface)
			@> \surfgrad >>
			H(\mathrm{curl},\thesurface)
			@> \surfcurl >>
			H(\mathrm{div},\thesurface)
			@> \surfdiv >>
			L^2(\thesurface)
			\\
			@VV{\pull_0}V @VV{\pull_1}V @VV{\pull_2}V @VV{\pull_3}V
			\\
			H^1(\thechart)
			@> \nabla >>
			H(\mathrm{curl},\thechart)
			@> \nabla \times >>
			H(\mathrm{div},\thechart)
			@> \nabla \cdot >>
			L^2(\thechart)
		\end{CD}
		\label{eq: de rham 3D}
	\end{equation}
	\begin{proof}
		Each identity follows from the definitions of \cref{sec: intrinsic operators}.
		The covariant components of the intrinsic gradient are $g \, ( g^{-1} \nabla f ) = \nabla f$.
		For the intrinsic curl, $\sqrt{g} \, \frac{1}{\sqrt{g}} \nabla \times \covar{v} = \nabla \times \covar{v}$.
		For the intrinsic divergence, $\sqrt{g} \, \frac{1}{\sqrt{g}} \nabla \cdot ( \sqrt{g} \, \contra{v} ) = \nabla \cdot ( \sqrt{g} \, \contra{v} ) = \nabla \cdot \flux{v}$. 
	\end{proof}
\end{proposition}

The parametric part of \cref{eq: de rham 3D} is the flat de Rham complex on $\thechart$, where the differential operators are the standard operators in flat Euclidean space.

\begin{proposition}
	\label{lem: duality}
	Let $\surfvec{v}, \surfvec{w}$ be tangent vector fields and $\surf{f}, \surf{\rho}$ scalar fields on $\thesurface$, with $n = 3$. Then
	\begin{align}
		\int_{\thesurface} \surfvec{v} \cdot \surfvec{w} \, \mathrm{d}\thesurface
		= \int_{\thechart} \pull_1 \surfvec{v} \cdot \pull_2 \surfvec{w} \, \mathrm{d}\chart{x} ,
		\qquad
		\int_{\thesurface} \surf{f} \, \surf{\rho} \, \mathrm{d}\thesurface
		= \int_{\thechart} \left( \pull_0 \surf{f} \right) \left( \pull_3 \surf{\rho} \right) \mathrm{d}\chart{x}.
		\label{eq: duality}
	\end{align}
	\begin{proof}
		By \cref{eq: inner product} and the definitions \cref{eq: pullbacks 3D}, pointwise on $\thechart$,
		$( \surfvec{v} \cdot \surfvec{w} ) \circ \map \, \sqrt{g} = \contra{v} \cdot g \, \contra{w} \, \sqrt{g} = \covar{v} \cdot ( \sqrt{g} \, \contra{w} ) = \pull_1 \surfvec{v} \cdot \pull_2 \surfvec{w}$
		and $( \surf{f} \, \surf{\rho} ) \circ \map \, \sqrt{g} = f  \sqrt{g} \rho  = ( \pull_0 \surf{f} ) ( \pull_3 \surf{\rho} )$. 
		Then \cref{eq: duality} follows from \cref{eq: integration}.
	\end{proof}
\end{proposition}

The duality pairing in \cref{lem: duality} is metric-free since the two arguments are represented in complementary {proxies}. That is, one argument is represented covariantly and the other as a flux or density.
In contrast, same-proxy pairings that have the same subscript carry the metric. That is, 
\begin{align}
	\int_{\thesurface} \surfvec{v} \cdot \surfvec{w} \, \mathrm{d}\thesurface
	= \int_{\thechart} \pull_1 \surfvec{v} \cdot g^{-1} \, \pull_1 \surfvec{w} \, \sqrt{g} \, \mathrm{d}\chart{x}
	= \int_{\thechart} \pull_2 \surfvec{v} \cdot g \, \pull_2 \surfvec{w} \, \frac{1}{\sqrt{g}} \, \mathrm{d}\chart{x} . 
	\label{eq: mass pairings}
\end{align}
These $g$-weighted products are precisely the mass terms of the \ac{fe} formulations of \cref{sec: models}.
Combining \cref{prop: commuting,lem: duality} with integration by parts on $\thechart$ yields the weak identities that underpin the compatible \discretizations\ of this work.
In what follows, $\int_{\thechart}$ abbreviates the sum of contributions over the partition domains of the atlas, and $\int_{\partial\thechart}$ the corresponding sum over their boundaries, where $\mathrm{d}s$ is the Euclidean boundary measure and $\covar n$ the Euclidean outward unit normal, regarded as a covariant field. 

\begin{proposition}
	\label{prop: ibp}
	Let $\surf{f} \in H^1(\thesurface)$ and $\surfvec{u} \in H(\mathrm{div},\thesurface)$, and set $f = \pull_0 \surf{f}$ and $\flux{u} = \pull_2 \surfvec{u}$. Then
	\begin{align}
		\int_{\thesurface} \surfgrad \surf{f} \cdot \surfvec{u} \, \mathrm{d}\thesurface
		+ \int_{\thesurface} \surf{f} \, \surfdiv \surfvec{u} \, \mathrm{d}\thesurface
		= \int_{\thechart} \nabla f \cdot \flux{u} \, \mathrm{d}\chart{x}
		+ \int_{\thechart} f \, \nabla \cdot \flux{u} \, \mathrm{d}\chart{x}
		= \int_{\partial \thechart} f \, \flux{u} \cdot \covar{n} \, \mathrm{d}s .
		\label{eq: ibp grad div}
	\end{align}
	When $n = 3$, for $\surfvec{v}, \surfvec{w} \in H(\mathrm{curl},\thesurface)$ with $\covar{v} = \pull_1 \surfvec{v}$ and $\covar{w} = \pull_1 \surfvec{w}$,
	\begin{equation}
	\begin{aligned}
		\int_{\thesurface} ( \surfcurl \surfvec{v} ) \cdot \surfvec{w} \, \mathrm{d}\thesurface
		- \int_{\thesurface} \surfvec{v} \cdot ( \surfcurl \surfvec{w} ) \, \mathrm{d}\thesurface
		&= \int_{\thechart} ( \nabla \times \covar{v} ) \cdot \covar{w} \, \mathrm{d}\chart{x}
		- \int_{\thechart} \covar{v} \cdot ( \nabla \times \covar{w} ) \, \mathrm{d}\chart{x} \\
		&= \int_{\partial \thechart} ( \covar{n} \times \covar{v} ) \cdot \covar{w} \, \mathrm{d}s .
		\label{eq: ibp curl}
	\end{aligned}
	\end{equation}
	\begin{proof}
		The ambient pairings in \cref{eq: ibp grad div,eq: ibp curl} reduce to Euclidean pairings on $\thechart$ by applying  \cref{lem: duality}  to the pairs $(\surfgrad \surf{f}, \surfvec{u})$, $(\surf{f}, \surfdiv \surfvec{u})$, $(\surfcurl \surfvec{v}, \surfvec{w})$, and $(\surfvec{v}, \surfcurl \surfvec{w})$, together with the commuting identities \cref{eq: commuting 3D}, where  the boundary terms are the standard Euclidean ones on $\thechart$.
	\end{proof}
\end{proposition}

\begin{remark}
	\label{rem: boundary terms}
	The boundary integrals in \cref{prop: ibp} involve only the Euclidean normal $\covar{n}$ to $\partial \thechart$. 
	When summing over the partition domains of a closed manifold, the interface contributions cancel and the boundary terms vanish altogether. 
	Boundary contributions only arise on genuine boundaries of $\thesurface$, such as the radial boundary of the \threeD atmospheric shell considered in \cref{sec: models} below.
\end{remark}

%%%%%%%%%%%%%%%%%%%%%%%%%%%%%%%%%%%%%%%%%%%%%%%%%%%%%%%%%%%%%%%%%%%%%%%%%%%%%%
\subsection{The rotated de Rham complex in two dimensions}
\label{sec: rotated complex}

Following the same pattern as the \threeD complex \cref{eq: de rham 3D}, we now define the \twoD rotated complex that is relevant to the shallow water equations in \cref{sec: models} below.
We collect here the corresponding definitions and identities.
For $n = 2$, we define the pullbacks
\begin{align}
	\pull_0 \surf{f} = f ,
	\qquad
	\pull_1 \surfvec{v} = \flux{v} ,
	\qquad
	\pull_2 \surf{\rho} = \sqrt{g} \, \rho ,
	\label{eq: pullbacks 2D}
\end{align}
where $\pull_1$ denotes the contravariant Piola map, so that in two dimensions the flux proxy $\flux{v} = \sqrt{g} \, \contra{v}$ carries the index one.
Sobolev spaces are defined chart-wise as in \cref{eq: conformity}, with $H(\mathrm{curl})$ omitted and $\pull_1$, $\pull_2$ indexed as in \cref{eq: pullbacks 2D}.

\begin{proposition}
	\label{prop: commuting 2D}
	For smooth fields and $n=2$,
	\begin{align}
		\pull_1 \left( \surfgrad^{\perp} \surf{f} \right) = \nabla^{\perp} \left( \pull_0 \surf{f} \right) ,
		\qquad
		\pull_2 \left( \surfdiv \surfvec{v} \right) = \nabla \cdot \left( \pull_1 \surfvec{v} \right) ,
		\label{eq: commuting 2D}
	\end{align}
such that the following diagram commutes:
	\begin{equation}
		\begin{CD}
			H^1(\thesurface)
			@> \surfgrad^{\perp} >>
			H(\mathrm{div},\thesurface)
			@> \surfdiv >>
			L^2(\thesurface)
			\\
			@VV{\pull_0}V @VV{\pull_1}V @VV{\pull_2}V
			\\
			H^1(\thechart)
			@> \nabla^{\perp} >>
			H(\mathrm{div},\thechart)
			@> \nabla \cdot >>
			L^2(\thechart)
		\end{CD}
		\label{eq: de rham 2D}
	\end{equation}
	\begin{proof}
		The proof follows as in \cref{prop: commuting}. For the intrinsic skew gradient, $\sqrt{g} \, \frac{1}{\sqrt{g}} \nabla^{\perp} f = \nabla^{\perp} f$.
	\end{proof}
\end{proposition}
 
{Under the indexing of \cref{eq: pullbacks 2D}, the scalar duality pairing in \cref{eq: duality} reads
$\int_{\thesurface} \surf{f} \, \surf{\rho} \, \mathrm{d}\thesurface = \int_{\thechart} ( \pull_0 \surf{f} ) ( \pull_2 \surf{\rho} ) \, \mathrm{d}\chart{x}$,
the pairing between two tangent fields is the flux form of \cref{eq: mass pairings}, namely
$\int_{\thesurface} \surfvec{v} \cdot \surfvec{w} \, \mathrm{d}\thesurface = \int_{\thechart} \flux{v} \cdot g \, \flux{w} \, \tfrac{1}{\sqrt{g}} \, \mathrm{d}\chart{x}$,
and \cref{eq: ibp grad div} holds verbatim with $\flux{u} = \pull_1 \surfvec{u}$.
The following skew integration by parts identity is genuinely \twoD.

\begin{proposition}
	\label{prop: ibp skew}
	For $\surf{\zeta} \in H^1(\thesurface)$ and a tangent field $\surfvec{u} \in [ L^2(\thesurface) ]^m$, with $\zeta = \pull_0 \surf{\zeta}$ and $\flux{u} = \pull_1 \surfvec{u}$,
	\begin{align}
		\int_{\thesurface} \surfvec{u} \cdot \surfgrad^{\perp} \surf{\zeta} \, \mathrm{d}\thesurface
		= - \int_{\thechart} \left( \sqrt{g} \, g^{-1} R \, \flux{u} \right) \cdot \nabla \zeta \, \mathrm{d}\chart{x} .
		\label{eq: ibp skew}
	\end{align}
	In addition, if $\surfgrad^{\perp} \cdot \surfvec{u} \in L^2(\thesurface)$, then
	\begin{align}
		{
		- \int_{\thechart} \left( \sqrt{g} \, g^{-1} R \, \flux{u} \right) \cdot \nabla \zeta \, \mathrm{d}\chart{x}
		= - \int_{\thesurface} \surf{\zeta} \, \surfgrad^{\perp} \cdot \surfvec{u} \, \mathrm{d}\thesurface
		- \int_{\partial \thechart} \zeta \left( \sqrt{g} \, g^{-1} R \, \flux{u} \right) \cdot \covar{n} \, \mathrm{d}s . }
		\label{eq: ibp skew boundary}
	\end{align}
	\begin{proof}
		The pairing is of same-proxy type \cref{eq: mass pairings}, so using $\pull_1 ( \surfgrad^{\perp} \surf{\zeta} ) = \nabla^{\perp} \zeta = R \nabla \zeta$ and $g R = \det(g) \, R \, g^{-1}$ \cite[Remark 15.11]{ErnGuermond_1}, yields
		\begin{align*}
			\int_{\thesurface} \surfvec{u} \cdot \surfgrad^{\perp} \surf{\zeta} \, \mathrm{d}\thesurface
			= \int_{\thechart} \flux{u} \cdot g \, ( R \nabla \zeta ) \, \frac{1}{\sqrt{g}} \, \mathrm{d}\chart{x}
			= - \int_{\thechart} \left( \sqrt{g} \, g^{-1} R \, \flux{u} \right) \cdot \nabla \zeta \, \mathrm{d}\chart{x} ,
		\end{align*}
		since $R^T = -R$.
		When $\surfgrad^{\perp} \cdot \surfvec{u} \in L^2(\thesurface)$,
		integrating by parts and noting that $\nabla \cdot ( \sqrt{g} \, g^{-1} R \, \flux{u} ) = - \pull_2 ( \surfgrad^{\perp} \cdot \surfvec{u} )$, by \cref{eq: skew ops}, yields \cref{eq: ibp skew boundary} via the second identity in \cref{eq: duality}.
	\end{proof}
\end{proposition}

\begin{remark}
	\label{rem: exact mimetics}
	
	The mimetic structure of the discrete de Rham complex is preserved exactly. 
	The compatible \discretizations\ of the model problems in \cref{sec: models}, built on the spaces of \cref{sec: fem}, require the commutation identities \cref{prop: commuting,prop: commuting 2D}, the duality pairings \cref{lem: duality}, the integration by parts identities \cref{prop: ibp}, and the annihilation identities \cref{prop: annihilation}. 
	These identities are all metric-free, and have piecewise polynomial integrands that can be computed exactly by numerical quadrature. 
	Furthermore, such \discretizations\ also require the $L^2$-type pairings \cref{eq: mass pairings} and the skew identity \cref{prop: ibp skew} for $n = 2$. 
	These metric-weighted pairings admit integration error for a non-polynomial metric, which reduces to machine precision as quantified in \cref{sec: quadrature study} below.
\end{remark}
\section{Intrinsic finite element spaces}
\label{sec: fem}

We now construct intrinsic finite element spaces.
In particular,
we discretise the parametric part of the commuting diagrams \cref{eq: de rham 3D,eq: de rham 2D}.  
This amounts to standard Euclidean \ac{fe} technology within each chart since 
inter-chart continuity is purely topological for all the spaces of the de Rham complex.
We subsequently discretise grad-conforming vector-valued spaces, which require \emph{transmission maps} to encode additional geometric data at chart interfaces.

%%%%%%%%%%%%%%%%%%%%%%%%%%%%%%%%%%%%%%%%%%%%%%%%%%%%%%%%%%%%%%%%%%%%%%%%%%%%%%
\subsection{Parametric triangulations and transmission maps}
\label{sec: triangulations}

Let $\mathcal{T}_k$ denote the \emph{chart triangulation} that is a conforming triangulation of the partition domain $P_k$ into cells obtained from a reference cell by affine maps.
For the cubed sphere manifold used to conduct experiments in \cref{sec: results} below, $\mathcal{T}_k$ is a Cartesian mesh of quadrilaterals ($n=2$) or hexahedra ($n=3$).
We denote  $\Gamma_{kl} = \map_k(\partial P_k) \cap \map_l(\partial P_l)$ as the interface between the image of two partition domains, which lies within the open overlap of the corresponding charts.
We assume the chart triangulations to be \emph{compatible} such that mesh entities (vertices, edges, faces) of $\mathcal{T}_k$ and $\mathcal{T}_l$ lying on $\Gamma_{kl}$ have matching images on $\thesurface$.
Then the collection $\mathcal{T} = \{ \mathcal{T}_k \}_{k=1}^{\kappa}$,  is the \emph{atlas triangulation}, together with the identification of the matching interface entities.
The image of $\mathcal{T}$ is a conforming mesh of $\thesurface$
by the compatibility of the chart triangulations assumed above.

In the multi-chart setting, the components of a tangent field have chart-wise definitions, which transform consistently across charts, as follows.
On the overlap of charts $k$ and $l$, the \emph{transmission map}, denoted $\transm_{kl}$, is the differential of the transition map $\tau_{kl}$. 
That is, 
\begin{align}
	\transm_{kl} = \mathrm{d}\tau_{kl} \in \mathbb{R}^{n \times n} ,
	\qquad
	\tau_{kl} = \map_l^{-1} \circ \map_k ,
	\qquad
	\contra{v}_l = \transm_{kl} \, \contra{v}_k ,
	\qquad
	\mathrm{d}\map_k = \mathrm{d}\map_l \, \transm_{kl} ,
	\label{eq: transmission}
\end{align}
where $\contra{v}_k$  and $\contra{v}_l$ denote the contravariant component fields of $\surfvec{v}$ in two different charts, the transmission map $\transm_{kl}$ represents a unique change-of-basis matrix such that $\contra{v}_l = \transm_{kl} \, \contra{v}_k $ when evaluated at corresponding points $\chart{x}$ and $\tau_{kl}(\chart{x})$, for which $\map_k(\chart{x}) = \map_l(\tau_{kl}(\chart{x}))$, and application of the chain rule yields the final expression in \cref{eq: transmission}. The transmission maps $\transm_{kl}$ are invertible since the transition maps $\tau_{kl}$ are diffeomorphisms on the overlaps.
The transmission maps satisfy the following properties at every interface point where they are defined:
\begin{enumerate}[label=(\roman*)]
	\item the \emph{cocycle conditions} $\transm_{kk} = I$ and $\transm_{lm} \transm_{kl} = \transm_{km}$, at points shared by more than two charts; 
	\item the \emph{metric compatibility} $g_k = \transm_{kl}^T \, g_l \, \transm_{kl}$;
	\item $\det \transm_{kl} > 0$, for a consistently oriented atlas.
\end{enumerate}
Property (ii) follows by substituting \cref{eq: transmission} into \cref{eq: metric}, and provides a useful consistency check between the transmission maps and the metrics of adjacent charts.

\begin{remark}
	\label{rem: transmission data}
	
	As shown in \cref{sec: grad conforming} below, the transmission maps are needed \emph{only} for the grad-conforming vector-valued spaces. The spaces of the discrete de Rham complexes in \cref{sec: discrete complex} require none.
	Since partition interfaces lie within open overlaps of the atlas, the transition maps are defined on neighbourhoods of the interfaces, and the transmission maps are fully determined by the atlas.  
	Equivalently, the transmission  maps can be computed through the embedding as $\transm_{kl} = (\mathrm{d}\map_l)^{+} \, \mathrm{d}\map_k$, where $(\cdot)^{+}$ denotes the Moore--Penrose pseudo-inverse.
	Resorting to the embedding is a convenience, not a requirement of the framework. 
\end{remark}

The following assumption concerns how the two charts parametrise a shared mesh entity.
\begin{assumption}
	\label{ass: entity parametrisation}
	The parametrisations of each shared entity induced by the affine cell maps of $\mathcal{T}_k$ and $\mathcal{T}_l$ agree up to a symmetry of the reference entity. Equivalently, the reparametrisation $\phi_l^{-1} \circ \tau_{kl} \circ \phi_k$ is affine, where $\phi_k$ and $\phi_l$ denote the \emph{entity parametrisations}. {These are restrictions of the affine cell maps to the corresponding entity of the reference cell.
	Thus, matching the vertices of the entity yields an element of its finite symmetry group. 
}
\end{assumption}
\cref{ass: entity parametrisation} is the standard conformity assumption of matched meshes in a single-domain \ac{fe} code, which is handled using typical orientation machinery \cite{Rognes2009,ErnGuermond_1}. 
For a vertex shared by two charts the condition in \cref{ass: entity parametrisation} is void, for a shared edge  it allows at most a reversal of direction, and for a shared facet a rotation or reflection.
In the multi-chart setting, \cref{ass: entity parametrisation} is a condition on the atlas that is satisfied whenever the transition maps are affine along the chart interfaces,
and is only needed for \ac{fe} spaces beyond the lowest order, as follows. 
It is also needed for the nodal construction of \cref{sec: grad conforming} for Lagrange spaces with order at least 2, where it ensures that both charts induce the same nodes on each shared entity.

%%%%%%%%%%%%%%%%%%%%%%%%%%%%%%%%%%%%%%%%%%%%%%%%%%%%%%%%%%%%%%%%%%%%%%%%%%%%%%
\subsection{Discrete de Rham complexes}
\label{sec: discrete complex}

On each chart triangulation $\mathcal{T}_k$ we consider standard \ac{fe} spaces that form a subcomplex of the Euclidean complex on $\thechart_k$. For $n = 3$, the subcomplex is
\begin{equation}
\resizebox{0.94\textwidth}{!}{$%
\begin{aligned}
	\V^0(\mathcal{T}_k) \subset H^1(\thechart_k)
	\; \xrightarrow{\; \nabla \;} \;
	\V^1(\mathcal{T}_k) \subset H(\mathrm{curl},\thechart_k)
	\; \xrightarrow{\; \nabla \times \;} \;
	\V^2(\mathcal{T}_k) \subset H(\mathrm{div},\thechart_k)
	\; \xrightarrow{\; \nabla \cdot \;} \;
	\V^3(\mathcal{T}_k) \subset L^2(\thechart_k) ,
	\label{eq: discrete complex}
\end{aligned}
$}
\end{equation}
Standard choices are continuous Lagrange, \N, \RT and discontinuous elements of compatible orders \cite{Nedelec1980,RT1977,Arnold2006}, in their tensor-product versions on quadrilateral and hexahedral meshes \cite{Arnold2005_quadrilateral,Arnold2018}. 
These spaces can be built from either full-polynomial or trimmed-polynomial families \cite{Arnold2006,Arnold2018}.
For the trimmed family indexed by $p$, the space $\V^0$ is the Lagrange space of order $p+1$, $\V^1$ is the \N space of order $p$, $\V^2$ is the \RT space of order $p$, and $\V^3$ is the discontinuous space of order $p$ \cite{Arnold2006}.
The lowest-order members of these families are the Whitney elements, where degrees of freedom attached to a mesh entity are
point values for $\V^0$, edge circulations $\int_e \covar{v} \cdot \mathrm{d}\chart{x}$ for $\V^1$, facet fluxes $\int_f \flux{u} \cdot \covar{n} \, \mathrm{d}s$ for $\V^2$, and cell integrals for $\V^3$, where $\covar{n}$ is the Euclidean unit normal to the facet $f$ in the parametric domain, as for $\partial\thechart$ in \cref{prop: ibp}. 
At arbitrary order, the entity degrees of freedom are \emph{moments} of the same quantities that are computed against polynomial weights defined through the entity parametrisation {of \cref{ass: entity parametrisation}}, together with interior moments, which are internal to a single chart and impose no inter-chart coupling.
The proceeding discussion applies for $n = 2$, where the rotated complex is $\V^0 \xrightarrow{\nabla^{\perp}} \V^1 \xrightarrow{\nabla \cdot} \V^2$ with continuous Lagrange, \RT and discontinuous elements.

Gluing the \emph{chart spaces $\V^s(\mathcal{T}_k)$} into \emph{global spaces} $\V^s(\mathcal{T})$ on the atlas triangulation requires the interface degrees of freedom from the two sides to agree.
Shared entities have dimension at most $n-1$. Cell degrees of freedom are interior to a single chart and require no identification.
The following proposition establishes the required \emph{pointwise} invariance on the shared entities, and involves no assumption on the parametrisations. Equality of the degrees of freedom, at any order, then follows by integration under \cref{ass: entity parametrisation}.
\begin{proposition}
	\label{prop: dof invariance}
	Consider two charts $k$ and $l$.
	Let $\surf{f}$, $\surfvec{v}$, $\surfvec{u}$ be fields on $\thesurface$, and denote   $f_k$, $\covar{v}_k$, $\flux{u}_k$ by their parametric representations \cref{eq: pullbacks 3D} in chart $k$ (indexed as in \cref{eq: pullbacks 2D} for $n=2$).
	Let $e_k$ and $e_l$  be corresponding interface entities of $\mathcal{T}_k$ and $\mathcal{T}_l$. That is a vertex, an edge ($n = 3$) or a facet under the interface identification. 
	Let $X_k$ be any parametrisation of $e_k$, so that $X_l = \tau_{kl} \circ X_k$ is a parametrisation of $e_l$. For corresponding facet parametrisations, let $J_{X_k} = \sqrt{\det( \mathrm{d}X_k^T \mathrm{d}X_k )}$, $J_{X_l} = \sqrt{\det( \mathrm{d}X_l^T \mathrm{d}X_l )}$, and let $\covar{n}_k$, $\covar{n}_l$ be the Euclidean unit normals to the facets $e_k$ and $e_l$ in the parametric domains of charts $k$ and $l$, oriented consistently under $\tau_{kl}$.  Then, pointwise in the parameter,
	\begin{align}
		f_k \circ X_k &= f_l \circ X_l ,
		&& \text{(vertices)}
		\nonumber \\
		( \covar{v}_k \circ X_k ) \cdot \mathrm{d}X_k &= ( \covar{v}_l \circ X_l ) \cdot \mathrm{d}X_l ,
		&& \text{(edges)}
		\nonumber \\
		\left( ( \flux{u}_k \cdot \covar{n}_k ) \circ X_k \right) J_{X_k} &= \left( ( \flux{u}_l \cdot \covar{n}_l ) \circ X_l \right) J_{X_l} .
		&& \text{(facets)}
		\label{eq: dof invariance}
	\end{align}
	\begin{proof}
		Let $\tau = \tau_{kl}$.
		Where the images of two charts overlap, $\map_k = \map_l \circ \tau$ implies, by the chain rule and \cref{eq: pullbacks 3D}, that the parametric representations are related by the transformations associated with $\tau$: 
		\begin{align}
			f_k = f_l \circ \tau ,
			\qquad
			\covar{v}_k = \mathrm{d}\tau^T ( \covar{v}_l \circ \tau ) ,
			\qquad
			\flux{u}_k = \det(\mathrm{d}\tau) \, \mathrm{d}\tau^{-1} ( \flux{u}_l \circ \tau ) ,
			\label{eq: functoriality}
		\end{align}
		These are  the identity, covariant and contravariant Piola transformations, respectively \cite[Chapter 9]{ErnGuermond_1}.
		The vertex identity in \cref{eq: dof invariance} is the first relation in \cref{eq: functoriality} evaluated along $X_k$.
		The edge identity in \cref{eq: dof invariance} is the second relation in \cref{eq: functoriality}, evaluated along $X_k$ and using $\tau \circ X_k = X_l$, which reads $\covar{v}_k \circ X_k = \mathrm{d}\tau^T ( \covar{v}_l \circ X_l )$, with $\mathrm{d}\tau$ evaluated along $X_k$. Hence
		\begin{align*}
			( \covar{v}_k \circ X_k ) \cdot \mathrm{d}X_k
			= \left( \mathrm{d}\tau^T ( \covar{v}_l \circ X_l ) \right) \cdot \mathrm{d}X_k
			= ( \covar{v}_l \circ X_l ) \cdot ( \mathrm{d}\tau \, \mathrm{d}X_k )
			= ( \covar{v}_l \circ X_l ) \cdot \mathrm{d}X_l ,
		\end{align*}
		where the second equality moves $\mathrm{d}\tau^T$ across the inner product, and the third is the chain rule $\mathrm{d}X_l = \mathrm{d}( \tau \circ X_k ) = \mathrm{d}\tau \, \mathrm{d}X_k$.
		The facet identity in \cref{eq: dof invariance} is the third relation in \cref{eq: functoriality}, evaluated along $X_k$ and using $\tau \circ X_k = X_l$, which reads $\flux{u}_k \circ X_k = \det( \mathrm{d}\tau ) \, \mathrm{d}\tau^{-1} ( \flux{u}_l \circ X_l )$, with $\mathrm{d}\tau$ evaluated along $X_k$. Writing the flux integrand as an inner product against the vector area element $\covar{n}_k J_{X_k}$,
		\begin{align*}
			\left( ( \flux{u}_k \cdot \covar{n}_k ) \circ X_k \right) J_{X_k}
			= ( \flux{u}_k \circ X_k ) \cdot ( \covar{n}_k \, J_{X_k} )
			&= ( \flux{u}_l \circ X_l ) \cdot \left( \det( \mathrm{d}\tau ) \, \mathrm{d}\tau^{-T} ( \covar{n}_k \, J_{X_k} ) \right)
			\\
			&= ( \flux{u}_l \circ X_l ) \cdot ( \covar{n}_l J_{X_l} )
			= \left( ( \flux{u}_l \cdot \covar{n}_l ) \circ X_l \right) J_{X_l} ,
		\end{align*}
		where the second equality moves $\mathrm{d}\tau^{-1}$ across the inner product, and the third is the Nanson identity $\det( \mathrm{d}\tau ) \, \mathrm{d}\tau^{-T} ( \covar{n}_k J_{X_k} ) = \covar{n}_l J_{X_l}$ that relates the vector area elements of corresponding facet parametrisations. For $n = 3$, the Nanson identity follows from applying $( \mathrm{d}\tau \, a ) \times ( \mathrm{d}\tau \, b ) = \det( \mathrm{d}\tau ) \, \mathrm{d}\tau^{-T} ( a \times b )$ to the columns of $\mathrm{d}X_k$, whose cross product is $\covar{n}_k J_{X_k}$.
	\end{proof}
\end{proposition}

\Cref{prop: dof invariance} is the key to gluing chart spaces into global spaces since it implies that the interface degrees of freedom in two different charts coincide, as we now make precise. Under \cref{ass: entity parametrisation}, take $X_k$ to be the entity parametrisation of $e_k$. Then $X_l$ is the entity parametrisation of $e_l$, up to a symmetry of the reference entity. Integrating \cref{eq: dof invariance} against the same polynomial weight then shows that the corresponding moments of the two charts coincide. In particular, the canonical degrees of freedom coincide at any order, with the common choice of weights per geometric entity realised exactly by the orientation machinery of a single-domain \ac{fe} code.
Identifying the matched interface degrees of freedom thus yields well-defined global spaces $\V^s(\mathcal{T})$, $s = 0, \dots, n$, at any order.
The fields on $\thesurface$ associated with the glued spaces are conforming in the sense of \cref{eq: conformity}, 
and form a discrete subcomplex, where the intrinsic derivative of a conforming field belongs to the next space of the complex. For example, $\V^0(\mathcal{T}) \subset H^1(\thesurface)$ whose derivatives belong to $\V^1(\mathcal{T}) \subset H(\mathrm{curl},\thesurface)$. 
In particular, since the glued spaces form a subcomplex and the annihilation identities of \cref{prop: annihilation} and the duality structure of \cref{lem: duality,prop: ibp} are algebraic identities of the Euclidean operators, which hold at the discrete level exactly, independently of the metric and of numerical quadrature, as anticipated in \cref{rem: exact mimetics}.
As in the standard single-domain \ac{fe} construction, the global field is built from cell-wise definitions that belong to the required space on each cell, and the shared degrees of freedom enforce the inter-cell trace continuity that characterises global conformity \cite{ErnGuermond_1}. 
Similarly at chart interfaces, the identification of the interface degrees of freedom  matches the traces between cells of adjacent charts.

\begin{remark}
	\label{rem: no transmission}
	
	The identities \cref{eq: dof invariance} are free of $\mathrm{d}\tau$, 
	meaning such identities do not require the metric, transmission maps or knowledge of the transition maps.
	In the edge identity, $\mathrm{d}\tau$ acts only along the entity, where $\mathrm{d}\tau \, \mathrm{d}X_k = \mathrm{d}X_l$ is determined by the two entity parametrisations, and in the facet identity the Piola factor $\det( \mathrm{d}\tau ) \, \mathrm{d}\tau^{-T}$ is exactly absorbed by the Nanson identity. In particular, the action of $\mathrm{d}\tau$ on directions transverse to the interface never enters.
\end{remark}

\begin{remark}
	\label{rem: lowest order} 
	For the lowest-order complex, \cref{ass: entity parametrisation} is not needed since the Whitney degrees of freedom are unweighted integrals, which are invariant under reparametrisation, and \cref{prop: dof invariance} applies to arbitrary parametrisations.
	The lowest-order complex therefore glues on any compatible atlas triangulation.
	Beyond the lowest order, \cref{ass: entity parametrisation} is required so that the two charts induce the same polynomial weights on each shared entity and the weighted moments coincide.
\end{remark}

\begin{remark}
	\label{rem: cohomology}
	The global discrete complex reproduces the de Rham cohomology of $\thesurface$ under the standard assumptions of \ac{fe} exterior calculus, namely the existence of bounded commuting projections \cite{Arnold2006,Arnold2018,HolstStern2012}.
	The specific manifolds used to perform numerical experiments in this study are a closed sphere when $n=2$ and a spherical shell when $n=3$.
	For the closed sphere, the cohomology is non-trivial only in $\V^0$ and $\V^2$, where the harmonic fields are the constant functions.
	The spherical shell is homotopy equivalent to the sphere, so its first cohomology group vanishes and the harmonic fields in $\V^0$ are again the constant functions.
	In particular, on both manifolds every tangent field with vanishing curl is a gradient.
	This determines the compatibility conditions on the source data of the Hodge Laplacian problems of \cref{sec: models}, which are analogous to the scalar and vector Poisson problems. On a closed manifold, the scalar Poisson problem is solvable provided $\int_{\thechart} f \sqrt{g} \, \mathrm{d}\chart{x} = 0$ holds, with the scalar unknown determined up to a constant. The discrete problems inherit the same constraints as their continuous counterparts.
\end{remark}

%%%%%%%%%%%%%%%%%%%%%%%%%%%%%%%%%%%%%%%%%%%%%%%%%%%%%%%%%%%%%%%%%%%%%%%%%%%%%%
\subsection{Grad-conforming vector-valued spaces}
\label{sec: grad conforming}

We now consider vector-valued Lagrangian chart spaces, $[\V^0(\mathcal{T}_k)]^n$, whose degrees of freedom are the \emph{pointwise values of the contravariant components} at the mesh nodes. Their gluing into a global space, denoted $\mathbb{W}(\mathcal{T})$, is the subject of this section.
In contrast with the entity integrals of \cref{prop: dof invariance} above, pointwise component values are \emph{not} chart-invariant. As per \cref{eq: transmission}, the values seen by two adjacent charts at an interface node $\surf{a}$ differ by the transmission map, 
\begin{align}
	\contra{v}_l(\chart{a}_l) = \transm_{kl}(\surf{a}) \, \contra{v}_k(\chart{a}_k) ,
	\label{eq: node transformation}
\end{align}
In contrast with \cref{rem: no transmission}, the action of the transmission map on directions transverse to the interface does not cancel, and is not determined by the entity correspondence and orientations.
Enforcing continuity of the associated tangent field therefore requires expressing the nodal degrees of freedom of all adjacent charts in one common frame, as illustrated in \cref{fig: FE}.
We now present the \emph{nodal construction} in \cref{sec: nodal construction} below for gluing vector-valued Lagrangian spaces across charts.
\begin{remark}
	\label{rem: stokes complex}
	Grad-conforming vector-valued spaces arise in complexes whose leading spaces demand $C^1$ elements of Hermite type. 
	These complexes are smoother than the de Rham complex in general.
	The discretisation of full smooth complexes, such as the Stokes complex \cite{Falk2013,Chen2024,Aznaran2022}, on manifolds is beyond the scope of this work. Here, grad-conforming vector fields are provided as a standalone capability.
\end{remark}

\begin{figure}[h!]
	\centering
	\includegraphics[width=0.95\textwidth]{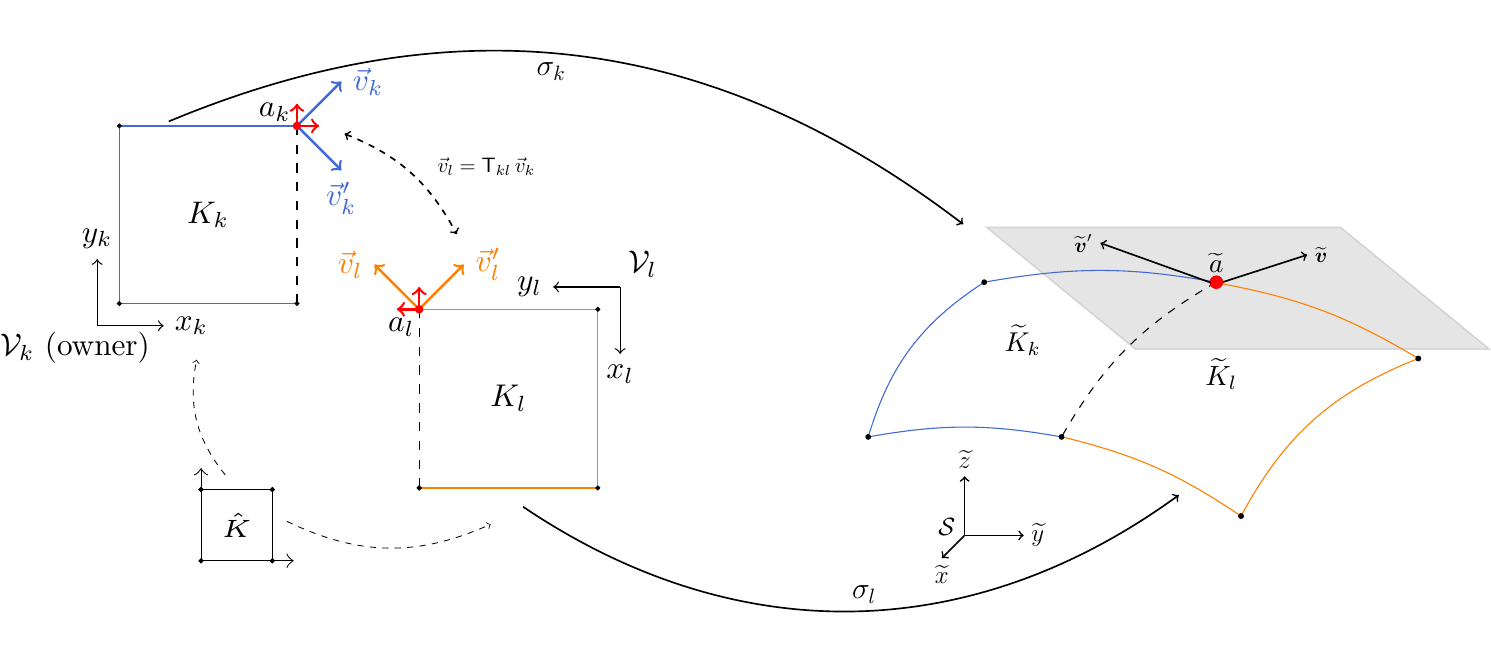}
	\caption{Two parametric cells in charts with different coordinate systems, mapped to ambient cells sharing an interface node $\surf{a}$ with parametric preimages $a_k$ and $a_l$.
	The nodal component values of tangent fields differ between charts by the transmission map, $\contra{v}_l = \transm_{kl} \contra{v}_k$ \cref{eq: node transformation}. The chart $\thechart_k$ is designated as the owner of $\surf{a}$.	}
	\label{fig: FE}
\end{figure}

\subsubsection{Nodal construction}
\label{sec: nodal construction}

For the nodal construction used in our experiments,
the glued grad-conforming vector-valued space $\mathbb{W}(\mathcal{T})$  consists of the fields $\surfvec{v}$ whose contravariant components belong to $[\V^0(\mathcal{T}_k)]^n$ on every chart, and whose nodal values satisfy \cref{eq: node transformation} at every interface node, which amounts to a non-standard gluing between chart spaces. 
The local space on every cell is the standard polynomial space of order $p+1$}$[\mathbb{Q}_{p+1}(\widehat{K})]^n$ \cite{Badia2018fempar}, for which the local basis, and hence the degrees of freedom, are modified through a change-of-basis matrix $C_K$, as follows. 

For each interface node $\surf{a}$, we designate an \emph{owner chart} $o(\surf{a})$ among the charts containing it, and define the global degrees of freedom at $\surf{a}$ as the components in the owner frame, $\contra{v}_{o(\surf{a})}(\chart{a})$.
For any other chart $k$ containing $\surf{a}$, the local nodal values are recovered as $\contra{v}_k = \transm_{o(\surf{a})k} \, \contra{v}_{o(\surf{a})}$ where the transmission map is evaluated at the node $\surf{a}$.
At the \ac{fe} level, this is a change of basis of the local shape functions. 
On a cell of chart $k$, the $n$ vector-valued shape functions attached to $\surf{a}$ are recombined by the change-of-basis block  $\transm_{o(\surf{a})k}$, so that their degrees of freedom become the owner components. The transpose of these blocks assemble into the cell change-of-basis matrix $C_K$, which is block-diagonal with identity blocks at interior and owner nodes. 
Then, the local basis
transforms as $\boldsymbol{\phi} \mapsto C_K \boldsymbol{\phi}$, while the element matrices and element vectors transform as  $A_K \mapsto C_K A_K C_K^T$ and $b_K \mapsto C_K b_K$ \cite{Badia2018fempar}.
Only cells touching a chart interface are affected, and the assembly is otherwise standard \cite{ErnGuermond_1}. 
The construction does not depend on the choice of owner since replacing $o(\surf{a})$ by another chart containing $\surf{a}$ is a change of basis of the global degrees of freedom at $\surf{a}$, which leaves the glued space unchanged.
Finally, at nodes shared by more than two charts, the pairwise relations imposed chart-wise  are mutually consistent by the cocycle condition (i) stated in \cref{sec: triangulations} above.
That is, $\transm_{kl} \transm_{o(\surf{a})k} = \transm_{o(\surf{a})l}$ for two non-owner charts $k$ and $l$, so the order of processing is irrelevant.

\begin{proposition}
	\label{prop: grad conforming}
	
	For the nodal construction, every field of $\mathbb{W}(\mathcal{T})$ defined above  is: (1) exactly tangent to $\thesurface$; (2) has components that are continuous, and thus $H^1$, within each chart; and (3) defines the same tangent vector from all adjacent charts at every interface node.
	\begin{proof}
		Tangency holds by construction, since the field on each chart collects the contravariant components from which the tangent field is generated through \cref{eq: contravariant expansion}. Continuity within a chart follows from the standard properties of Lagrangian elements.
		At an interface node, \cref{eq: node transformation} and \cref{eq: transmission} imply $\mathrm{d}\map_k \, \contra{v}_k = \mathrm{d}\map_l \, \contra{v}_l$. That is, the two charts define the same tangent vector.
		This holds at all shared nodes by the compatibility of the atlas triangulation together with \cref{ass: entity parametrisation}, which makes the two entity parametrisations agree up to a symmetry of the reference entity.
		Since the Lagrange node set is invariant under that symmetry group, both charts induce the same nodes on each shared entity, up to the permutation of the local numbering handled by the usual orientation machinery.
		For $p = 0$ the nodes are the entity vertices and the compatibility of the atlas triangulation suffices.
	\end{proof}
\end{proposition}

\begin{corollary}
	\label{cor: constant transmission} 
	If the transmission maps are constant along each interface, the fields of $\mathbb{W}(\mathcal{T})$ are single-valued on entire shared entities and $\mathbb{W}(\mathcal{T})$ is grad-conforming in the sense of \cref{eq: conformity} applied componentwise. This holds, in particular, for single-chart manifolds, where no interfaces are present, and for atlases whose transition maps are affine in a neighbourhood of the interfaces.
	\begin{proof}
		Along a shared entity, the traces of $\contra{v}_k$ and $\contra{v}_l$ are polynomial in their own frames, and single-valuedness of the tangent field amounts to the relation $\contra{v}_l = \transm_{kl} \, \contra{v}_k$ holding \emph{pointwise} along the entity, whereas the nodal identification enforces the relation only at the nodes, where the transmission map is evaluated.
		When $\transm_{kl}$ is constant along the entity, $\transm_{kl} \, \contra{v}_k$ is itself polynomial and agrees with $\contra{v}_l$ at the entity nodes, hence everywhere on the entity, by the nodal determination of Lagrangian elements.
		Single-valuedness on the interfaces, together with chart-wise continuity, yields \cref{eq: conformity} applied componentwise.
		For a single-chart manifold there are no interfaces and the statement is immediate. For atlases with affine transition maps near the interfaces, the transmission maps are constant there.
	\end{proof}
\end{corollary}

\begin{remark}
	\label{rem: tangent H1}
	Within each chart, and whenever the transmission maps are constant along the interfaces (\cref{cor: constant transmission}), the glued space is simultaneously exactly tangent to $\thesurface$ and $H^1$-conforming since the geometry is exact and smooth.
	When the transmission maps vary along an interface, single-valuedness between the nodes fails in general and the glued space is not fully $H^1$-conforming. 
	This is the case for the equiangular cubed sphere manifold outlined \cref{app: cubed sphere}, where the transmission map is not constant along the interface of charts. That is, on the edge shared by two adjacent panels, with $x^2$ the coordinate along the edge, the panel maps stated in \citet{Nair2005} yield $\transm_{kl} = \left[ \begin{smallmatrix} 1 & 0 \\ -\sin 2 x^2 & 1 \end{smallmatrix} \right]$, which is not constant along the edge.
	\cref{ass: entity parametrisation}, which the cubed sphere atlas satisfies, does not prevent this failure. The assumption constrains only the transition maps \emph{restricted to} the interfaces, whereas the transmission maps involve the derivatives of the transition maps in the direction transverse to the interface, which \cref{ass: entity parametrisation} leaves unconstrained. 
\end{remark}

\begin{remark}
	\label{rem: accuracy}
	On a shared entity $e$, the two traces are $\contra{w}_l$ and $\transm_{kl} \contra{w}_k$, which agree at the entity nodes, so that $\contra{w}_l = I_{p+1} ( \transm_{kl} \contra{w}_k )$ on $e$, where $I_{p+1}$ denotes the nodal interpolation of degree $p+1$ on the entity.
	Since $I_{p+1}$ reproduces polynomials of degree $p+1$,
	\begin{align*}
		\transm_{kl} \contra{w}_k - \contra{w}_l
		= ( I - I_{p+1} ) ( \transm_{kl} \contra{w}_k )
		= ( I - I_{p+1} ) ( \transm_{kl} \contra{w}_k - \pi )
		\qquad \text{on } e ,
	\end{align*}
	for every polynomial $\pi$ of degree $p+1$ on $e$.
	Taking the infimum over $\pi$ and applying the Bramble--Hilbert lemma yields
	\begin{align}
		h_e^{-1/2} \, \| \transm_{kl} \contra{w}_k - \contra{w}_l \|_{L^2(e)}
		\leq C \, h_e^{\, p + 3/2} \sum_{m=0}^{p+1}
		\| \partial_s^{\, p+2-m} \transm_{kl} \|_{L^{\infty}(e)} \,
		\| \partial_s^{\, m} \contra{w}_k \|_{L^2(e)} ,
		\label{eq: jump bound}
	\end{align}
	where $s$ is a coordinate along the entity.
\end{remark}

%%%%%%%%%%%%%%%%%%%%%%%%%%%%%%%%%%%%%%%%%%%%%%%%%%%%%%%%%%%%%%%%%%%%%%%%%%%%%%
\subsection{Assembly}
\label{sec: assembly}

All element integrals are computed on the Cartesian meshes of the parametric domains. 
Each cell map is affine with a constant Jacobian, and the metric quantities $g$, $g^{-1}$, $\sqrt{g}$ are evaluated at the quadrature points of the reference cell. 
This computation is analytical for manifold atlases in closed form.  
By \cref{rem: exact mimetics}, quadrature accuracy affects only the metric-weighted integrals, not the mimetic structure of the discrete complex.
No object of the ambient dimension $m$ is ever assembled.
The computational implications of this intrinsic pipeline, in comparison with extrinsic assembly, are quantified in \cref{sec: results} below.

\section{Intrinsic formulations of model problems}
\label{sec: models}

We now derive the intrinsic \ac{fe} formulations for the Hodge Laplacian problems associated to the de Rham complex, and the rotating shallow water equations.
Throughout, $\V^s(\mathcal{T})$ denotes the glued global spaces defined in \cref{sec: discrete complex} above, with chart-interface terms cancelling via \cref{rem: boundary terms}.
Problems are posed at the discrete level, 
where replacing discrete spaces by their Sobolev counterparts \cref{eq: conformity} yields corresponding continuous weak forms.
We subsequently consider the surface Stokes problem that involves grad-conforming vector-valued spaces developed in \cref{sec: grad conforming}.

%%%%%%%%%%%%%%%%%%%%%%%%%%%%%%%%%%%%%%%%%%%%%%%%%%%%%%%%%%%%%%%%%%%%%%%%%%%%%%
\subsection{Hodge Laplacian problems}
\label{sec: hodge}

The de Rham complex carries a family of abstract Hodge Laplacian problems \cite[Chapter 4]{Arnold2018}.
For $n = 3$, we consider the mixed scalar Poisson and vector Poisson problems, 
posed on the $H(\mathrm{div})$--$L^2$ and $H^1$--$H(\mathrm{curl})$ pair of \cref{eq: de rham 3D}, respectively.

\subsubsection{Scalar Poisson problem}
\label{sec: scalar poisson}

Given $\surf{f} : \thesurface \rightarrow \mathbb{R}$, the mixed Poisson problem for the Laplace--Beltrami operator on $\thesurface$ with Dirichlet boundary conditions, is: find $\surfvec{u}$ tangent and $\surf{\varphi} : \thesurface \rightarrow \mathbb{R}$ such that
\begin{align}
	\surfvec{u} + \surfgrad \surf{\varphi} = \boldsymbol{0} ,
	\qquad
	\surfdiv \surfvec{u} = \surf{f}
	\quad \text{on } \thesurface ,
	\qquad
	\surf{\varphi} = \surf{c}
	\quad \text{on } \partial\thesurface .
	\label{eq: poisson mixed}
\end{align}
Representing the flux as $\flux{u} = \pull_2 \surfvec{u}$ and the scalar unknown by its pullback $\varphi = \pull_0 \surf{\varphi}$, the intrinsic \ac{fe} formulation of \cref{eq: poisson mixed} is:
find $\flux{u}_h \in \V^2(\mathcal{T})$, $\varphi_h \in \V^3(\mathcal{T})$ such that
\begin{subequations}
	\label{eq: poisson weak}
	\begin{align}
		\int_{\thechart} \flux{u}_h \cdot g \, \flux{\psi} \, \frac{1}{\sqrt{g}} \, \mathrm{d}\chart{x}
		- \int_{\thechart} \varphi_h \, \nabla \cdot \flux{\psi} \, \mathrm{d}\chart{x}
		&= - \int_{\partial\thechart} c \, \flux{\psi} \cdot \covar{n} \, \mathrm{d}s ,
		&& \forall \flux{\psi} \in \V^2(\mathcal{T}) ,
		\label{eq: poisson weak u}
		\\
		\int_{\thechart} \omega \, \nabla \cdot \flux{u}_h \, \mathrm{d}\chart{x}
		&= \int_{\thechart} f \, \omega \, \sqrt{g} \, \mathrm{d}\chart{x} ,
		&& \forall \omega \in \V^3(\mathcal{T}) .
		\label{eq: poisson weak phi}
	\end{align}
\end{subequations}
This formulation is a direct application of the machinery of \cref{sec: pullbacks} above. 
In \cref{eq: poisson weak u},
the first term is the flux mass pairing \cref{eq: mass pairings}, while the second term arises from applying the duality pairing \cref{eq: duality} to $(\surfgrad\surf{\varphi}, \surfvec{\psi})$, the commutation \cref{eq: commuting 3D},  and the Euclidean integration by parts \cref{eq: ibp grad div}. 
\cref{eq: poisson weak phi} follows from \cref{eq: duality,eq: commuting 3D}.
The metric is confined to the mass term, while the divergence terms are metric-free.

\begin{remark}
	\label{rem: L2 representation}
	In \cref{eq: poisson weak} the  space $\V^3$ describes the plain pullback $\varphi_h$, not the density $\pull_3 \surf{\varphi} = \sqrt{g}\,\varphi$.
	Both choices yield conforming discretisations, by the boundedness of the metric weights, but they span different discrete spaces. The former is adopted hereinafter.
\end{remark}

\subsubsection{Vector Poisson problem}
\label{sec: vector poisson}

{In mixed form, the vector Poisson problem reads:} find $\surfvec{q}$ tangent and $\surf{w} : \thesurface \rightarrow \mathbb{R}$ such that  
\begin{align}
	\surf{w} + \surfdiv \surfvec{q} = 0 ,
	\qquad
	\surfcurl ( \surfcurl \surfvec{q} ) + \surfgrad \surf{w} = \surfvec{f}
	\quad \text{on } \thesurface ,
	\label{eq: vector poisson}
\end{align}
subject to prescribed values of the normal trace of $\surfvec{q}$ and of the tangential trace of $\surfcurl \surfvec{q}$ on $\partial\thesurface$, made precise in \cref{rem: boundary data} below.
Representing $\covar{q} = \pull_1 \surfvec{q}$ and $w = \pull_0 \surf{w}$, and assuming $\thesurface$ has no tunnels, so that every tangent field with vanishing surface curl is a surface gradient  (\cref{rem: cohomology}), the intrinsic formulation is:
find $w_h \in \V^0(\mathcal{T})$, $\covar{q}_h \in \V^1(\mathcal{T})$ such that
\begin{subequations}
	\label{eq: vector poisson weak}
	\begin{align}
		\int_{\thechart} w_h \, \xi \, \sqrt{g} \, \mathrm{d}\chart{x}
		&- \int_{\thechart} \nabla \xi \cdot g^{-1} \covar{q}_h \, \sqrt{g} \, \mathrm{d}\chart{x}
		= - \int_{\partial\thechart} \xi \, \chi \, \mathrm{d}s ,
		&& \forall \xi \in \V^0(\mathcal{T}) ,
		\label{eq: vector poisson weak sigma}
		\\
		\int_{\thechart} ( \nabla \times \covar{q}_h ) \cdot g \, ( \nabla \times \covar{\omega} ) \, \frac{1}{\sqrt{g}} \, \mathrm{d}\chart{x}
		&+ \int_{\thechart} \nabla w_h \cdot g^{-1} \covar{\omega} \, \sqrt{g} \, \mathrm{d}\chart{x}
		\nonumber \\
		&= \int_{\thechart} \covar{f} \cdot g^{-1} \covar{\omega} \, \sqrt{g} \, \mathrm{d}\chart{x}
		- \int_{\partial\thechart} ( \covar{n} \times \covar{\chi} ) \cdot \covar{\omega} \, \mathrm{d}s ,
		&& \forall \covar{\omega} \in \V^1(\mathcal{T}) ,
		\label{eq: vector poisson weak q}
	\end{align}
\end{subequations}
where $\covar{f} = \pull_1 \surfvec{f}$, and $\covar{n}$ is the Euclidean normal to $\partial \thechart$.
The mass term in \cref{eq: vector poisson weak sigma} and the gradient term $\nabla w_h \cdot g^{-1} \covar{\omega}$ in \cref{eq: vector poisson weak q} are covariant same-proxy pairings \cref{eq: mass pairings}. The term $\nabla\xi \cdot g^{-1}\covar{q}_h$ in \cref{eq: vector poisson weak sigma} arises from the duality pairing of $(\surfgrad\surf{\xi}, \surfvec{q})$ after integration by parts \cref{eq: ibp grad div}. 
The curl--curl term in 	\cref{eq: vector poisson weak q} is the flux same-proxy pairing of $\nabla\times\covar{q}_h$, which is obtained from the curl identity \cref{eq: ibp curl} and the commutation \cref{eq: commuting 3D}. 

\begin{remark}
	\label{rem: boundary data}
	Boundary data enters \cref{eq: vector poisson weak} as \emph{parametric densities} where $\chi$ is the parametric representation of the normal trace datum  of $\surfvec{q}$ per unit parametric boundary length, and $\covar{\chi}$ is the covariant representation of the tangential trace  of $\surfcurl\surfvec{q}$.
	Converting ambient data into $\chi$ and $\covar{\chi}$ is a data-preparation step that is available in closed form
	 for the manufactured solutions in \cref{sec: results} below.
	The same observation applies to the Dirichlet datum $c$ in \cref{eq: poisson weak u}.
\end{remark}

%%%%%%%%%%%%%%%%%%%%%%%%%%%%%%%%%%%%%%%%%%%%%%%%%%%%%%%%%%%%%%%%%%%%%%%%%%%%%%
\subsection{Rotating shallow water equations}
\label{sec: swe}

The rotating shallow water equations on a closed \twoD manifold $\thesurface$ describe the evolution of a tangent velocity field $\surfvec{u}$ and of the fluid depth $\surf{\varphi} : \thesurface \rightarrow \mathbb{R}$ as
\begin{subequations}
	\label{eq: swe}
	\begin{align}
		\partial_t \surfvec{u}
		+ \surf{q} \, \surfvec{F}^{\perp}
		+ \surfgrad \surf{\Phi} &= \boldsymbol{0}
		\qquad \text{in } \thesurface \times (0, t] ,
		\label{eq: swe u}
		\\
		\partial_t \surf{\varphi}
		+ \surfdiv \surfvec{F} &= 0
		\qquad \text{in } \thesurface \times (0, t] ,
		\label{eq: swe phi}
		\\
		\text{with diagnostic variables}\qquad 
		\surfvec{F} = \surf{\varphi} \, \surfvec{u} ,
		\qquad
		\surf{\Phi} &= \tfrac{1}{2} \, \surfvec{u} \cdot \surfvec{u} + \gravity \, ( \surf{\varphi} + \surf{\topography} ) ,
		\qquad
		\surf{q} = \frac{ \surfgrad^{\perp} \cdot \, \surfvec{u} + \surf{\coriolis} }{ \surf{\varphi} } ,
		\label{eq: swe diagnostics}
	\end{align}
\end{subequations}
which are the mass flux, the Bernoulli potential and the potential vorticity, respectively. 
In \cref{eq: swe diagnostics}, $\surf{\topography}$ is the topography, $\gravity$ the acceleration of gravity and $\surf{\coriolis}$ the Coriolis parameter.
Here $\surfvec{v}^{\perp}$ denotes the pointwise rotation of a tangent field by $\pi/2$. 
Such a rotation is performed consistently with the orientation of $\thesurface$, and is defined intrinsically as follows.

\begin{proposition}
	\label{lem: rotation}
	Let $n = 2$ and let $\surfvec{v}$ be a tangent field with flux proxy $\flux{v} = \pull_1 \surfvec{v}$. The rotated field $\surfvec{v}^{\perp}$ has flux proxy
	\begin{align}
		\pull_1 ( \surfvec{v}^{\perp} ) = \sqrt{g} \, g^{-1} R \, \flux{v} ,
		\label{eq: rotation proxy}
	\end{align}
	and, for any tangent fields $\surfvec{v}, \surfvec{w}$ and scalar $\surf{a}$, the rotated pairing is metric-free:
	\begin{align}
		\int_{\thesurface} \surf{a} \; \surfvec{v}^{\perp} \cdot \, \surfvec{w} \, \mathrm{d}\thesurface
		= \int_{\thechart} a \, ( R \flux{v} ) \cdot \flux{w} \, \mathrm{d}\chart{x} .
		\label{eq: rotation pairing}
	\end{align}
	\begin{proof}
		
		Let us consider the tangent field $\surfvec{z}$ whose contravariant components are $\contra{z} = \sqrt{g} \, g^{-1} R \, \contra{v}$, where $\contra{v} = \flux{v}/\sqrt{g}$ are the contravariant components of $\surfvec{v}$. 
		Recall $R^T = -R$, $R^T g^{-1} R = g / \det g$, and \cref{eq: inner product}. 
		Then, the field $\surfvec{z}$  is pointwise orthogonal to $\surfvec{v}$, of equal norm and consistently oriented. 	
		Multiplying  the contravariant components $\contra{z}$ by $\sqrt{g}$ gives the flux proxy of $\surfvec{z}$, which is \cref{eq: rotation proxy}.
		For \cref{eq: rotation pairing}, the same-proxy pairing \cref{eq: mass pairings} yields
		$\int_{\thechart} a \, ( \sqrt{g} g^{-1} R \flux{v} ) \cdot g \flux{w} \, \tfrac{1}{\sqrt{g}} \mathrm{d}\chart{x}
		= \int_{\thechart} a \, ( R \flux{v} ) \cdot \flux{w} \, \mathrm{d}\chart{x}$.
	\end{proof}
\end{proposition}

\subsubsection{Semi-discrete formulation}
\label{sec: swe semidiscrete}

Let $\flux{u} = \pull_1 \surfvec{u}$, $\flux{F} = \pull_1 \surfvec{F}$, and let $\varphi, \Phi, q, \topography, \coriolis$ denote plain pullbacks as in \cref{rem: L2 representation}.
Using \cref{eq: de rham 2D}, the intrinsic \ac{fe} formulation of the prognostic variables is: find  $\flux{u}_h(t) \in \V^1(\mathcal{T})$, $\varphi_h(t) \in \V^2(\mathcal{T})$ such that
\begin{subequations}
	\label{eq: swe weak}
	\begin{align}
		\int_{\thechart} \partial_t \flux{u}_h \cdot g \, \flux{\psi} \, \frac{1}{\sqrt{g}} \, \mathrm{d}\chart{x}
		+ \int_{\thechart} q^{\mathrm{u}}_h \, ( R \flux{F}_h ) \cdot \flux{\psi} \, \mathrm{d}\chart{x}
		- \int_{\thechart} \Phi_h \, \nabla \cdot \flux{\psi} \, \mathrm{d}\chart{x}
		&= 0 ,
		&& \forall \flux{\psi} \in \V^1(\mathcal{T}) ,
		\label{eq: swe weak u}
		\\
		\int_{\thechart} \partial_t \varphi_h \, \omega \, \sqrt{g} \, \mathrm{d}\chart{x}
		+ \int_{\thechart} \omega \, \nabla \cdot \flux{F}_h \, \mathrm{d}\chart{x}
		&= 0 ,
		&& \forall \omega \in \V^2(\mathcal{T}) .
		\label{eq: swe weak phi}
	\end{align}
\end{subequations}
The potential vorticity is upwinded via the streamline upwind Petrov--Galerkin method \cite{BrooksHughes1982,Lee2024},
\begin{align}
	q^{\mathrm{u}}_h = q_h - \vartheta \left( \partial_t q_h + \frac{\flux{u}_h}{\sqrt{g}} \cdot \nabla q_h \right) ,
	\label{eq: supg}
\end{align}
where $\vartheta$ is the stabilisation parameter, and the advective term in \cref{eq: supg} is the pullback of $\surfvec{u} \cdot \surfgrad \surf{q}$.
In \cref{eq: swe weak u},
the Coriolis term  is metric-free by \cref{lem: rotation}, and the Bernoulli term follows from \cref{eq: ibp grad div} on the closed manifold.
\cref{eq: swe weak phi} follows from \cref{eq: duality,eq: commuting 2D}.

The diagnostic variables $\flux{F}_h \in \V^1(\mathcal{T})$, $\Phi_h \in \V^2(\mathcal{T})$, $q_h \in \V^0(\mathcal{T})$ solve
\begin{subequations}
	\label{eq: swe diagnostics weak}
	\begin{align}
		\int_{\thechart} \flux{F}_h \cdot g \, \flux{\psi} \, \frac{1}{\sqrt{g}} \, \mathrm{d}\chart{x}
		&= \int_{\thechart} \varphi_h \, \flux{u}_h \cdot g \, \flux{\psi} \, \frac{1}{\sqrt{g}} \, \mathrm{d}\chart{x} ,
		&& \forall \flux{\psi} \in \V^1(\mathcal{T}) ,
		\label{eq: weak F}
		\\
		\int_{\thechart} \Phi_h \, \omega \, \sqrt{g} \, \mathrm{d}\chart{x}
		&= \int_{\thechart} \frac{1}{2} \, \flux{u}_h \cdot g \, \flux{u}_h \, \omega \, \frac{1}{\sqrt{g}} \, \mathrm{d}\chart{x}
		+ \int_{\thechart} \gravity \, ( \varphi_h + \topography ) \, \omega \, \sqrt{g} \, \mathrm{d}\chart{x} ,
		&& \forall \omega \in \V^2(\mathcal{T}) ,
		\label{eq: weak Phi}
		\\
		\int_{\thechart} q_h \, \varphi_h \, \xi \, \sqrt{g} \, \mathrm{d}\chart{x}
		&= \int_{\thechart} \left( \sqrt{g} \, g^{-1} R \, \flux{u}_h \right) \cdot \nabla \xi \, \mathrm{d}\chart{x}
		+ \int_{\thechart} \coriolis \, \xi \, \sqrt{g} \, \mathrm{d}\chart{x} ,
		&& \forall \xi \in \V^0(\mathcal{T}) ,
		\label{eq: weak q}
	\end{align}
\end{subequations}
where the first term on the right of \cref{eq: weak q} is the weak form of the relative vorticity $\surfgrad^{\perp} \cdot \, \surfvec{u}$, which is obtained via the skew identity \cref{eq: ibp skew} on the closed manifold.

\subsubsection{Conservation properties}
\label{sec: swe conservation}

The semi-discrete formulation \cref{eq: swe weak,eq: swe diagnostics weak} maintains the skew-symmetric structure of the  continuous equations \cite{Cotter2012,Wimmer2020,BauerCotter2018}, and yields the following conservation properties.

\begin{proposition}
	\label{prop: mass conservation}
	The semi-discrete system \cref{eq: swe weak} conserves the discrete mass $\mathcal{M}_h(t) = \int_{\thechart} \varphi_h \sqrt{g} \, \mathrm{d}\chart{x}$.
	\begin{proof}
		Take $\omega = 1$ in \cref{eq: swe weak phi}. 
		Applying \cref{rem: boundary terms,prop: dof invariance} yields
		 $\mathrm{d}\mathcal{M}_h / \mathrm{d}t = - \int_{\thechart} \nabla \cdot \flux{F}_h \, \mathrm{d}\chart{x} = 0$ 
		 since $\flux{F}_h\in \mathbb{V}^1(\mathcal{T})$.
	\end{proof}
\end{proposition}

\begin{proposition}
	\label{prop: energy conservation}
	The semi-discrete system \cref{eq: swe weak,eq: swe diagnostics weak} conserves the discrete energy
	\begin{align}
		\mathcal{H}_h(t) =
		\int_{\thechart} \frac{1}{2} \, \varphi_h \, \flux{u}_h \cdot g \, \flux{u}_h \, \frac{1}{\sqrt{g}} \, \mathrm{d}\chart{x}
		+ \int_{\thechart} \gravity \left( \frac{1}{2} \, \varphi_h^2 + \varphi_h \topography \right) \sqrt{g} \, \mathrm{d}\chart{x} .
		\label{eq: energy}
	\end{align}
	\begin{proof}
		Differentiating \cref{eq: energy} in time and using \cref{eq: weak F,eq: weak Phi} yields
		\begin{align*}
			\frac{\mathrm{d}\mathcal{H}_h}{\mathrm{d}t}
			= \int_{\thechart} \partial_t \flux{u}_h \cdot g \, \flux{F}_h \, \frac{1}{\sqrt{g}} \, \mathrm{d}\chart{x}
			+ \int_{\thechart} \partial_t \varphi_h \, \Phi_h \, \sqrt{g} \, \mathrm{d}\chart{x} ,
		\end{align*}
		since the functional derivatives of $\mathcal{H}_h$ with respect to $\flux{u}_h$ and $\varphi_h$ are represented by $\flux{F}_h$ and $\Phi_h$.
		Taking $\flux{\psi} = \flux{F}_h$ in \cref{eq: swe weak u} and $\omega = \Phi_h$ in \cref{eq: swe weak phi} yields $\mathrm{d}\mathcal{H}_h / \mathrm{d}t = 0$
		since the Coriolis contribution in \cref{eq: swe weak u} vanishes pointwise for any coefficient $q^{\mathrm{u}}_h$	by the skew-symmetry of $R$, and the two Bernoulli contributions cancel.
	\end{proof}
\end{proposition}

\begin{remark}
	\label{rem: energy topography}
	\Cref{prop: energy conservation} holds with time-independent topography included. For $\topography = 0$, \cref{eq: energy} reduces to the standard discrete energy that is the sum of kinetic and potential energy. 
	Both conservation statements are consequences of the metric-free structure highlighted in \cref{rem: exact mimetics}. Consequently, mass and energy conservation hold exactly regardless of the accuracy of the metric quadrature. 
\end{remark}

\subsubsection{Extension to three dimensions}
\label{sec: swe 3d}

We now let $\thesurface$ be a \threeD atmospheric shell of co-dimension zero that is discretised with a single element in the radial direction. 
In this case, the shallow water equations can be expressed as
\begin{subequations}
	\label{eq: swe 3d}
	\begin{align}
		\partial_t \surfvec{u}
		+ \surfvec{q} \times \surfvec{F}
		+ \surfgrad \surf{\Phi} &= \boldsymbol{0}
		\qquad \text{in } \thesurface \times (0, t] ,
		\\
		\partial_t \surf{\varphi}
		+ \surfdiv \surfvec{F} &= 0
		\qquad \text{in } \thesurface \times (0, t] ,
		\\
		\text{with diagnostic variables} \quad
		\surfvec{F} = \surf{\varphi} \, \surfvec{u} ,
		\qquad
		\surf{\Phi} = \tfrac{1}{2} \, \surfvec{u} \cdot \surfvec{u} &+ \gravity \, ( \surf{\varphi} + \surf{\topography} ) ,
		\qquad
		\surfvec{q} = \frac{ \surfcurl \surfvec{u} + \surf{\coriolis} \, \surfvec{k} }{ \surf{\varphi} } ,
		\label{eq: swe 3d diagnostics}
	\end{align}
\end{subequations}
Here $\surfvec{k}$ is the outward unit vector in the radial direction, so that $\surf{\coriolis} \, \surfvec{k}$ is the radial Coriolis field whose covariant representation appears in \cref{eq: swe 3d coriolis} below.
In parametric form,
the conditions on the radial boundaries are $\flux{u} \cdot \covar{n} = 0$ and $\covar{n} \times \covar{q} = \contra{0}$ on $\partial\thechart$.
The cross products in \cref{eq: swe 3d} have the following intrinsic representation.

\begin{proposition}
	\label{lem: cross product} 
	Let $n = m = 3$ and let $\surfvec{v}, \surfvec{w}$ be fields with covariant representations $\covar{v}, \covar{w}$. Then $\surfvec{v} \times \surfvec{w}$ has contravariant components $( \covar{v} \times \covar{w} ) / \sqrt{g}$.
	\begin{proof}
		With $A = \mathrm{d}\map$ invertible and $\det A = \sqrt{g}$ for a positively oriented chart, the identity $(A\contra{a}) \times (A\contra{b}) = \det(A) \, A^{-T} ( \contra{a} \times \contra{b} )$ and $\covar{v} = g \contra{v}$ give
		$\surfvec{v} \times \surfvec{w} = \sqrt{g} \, A^{-T} ( \contra{v} \times \contra{w} )$ and
		$\covar{v} \times \covar{w} = \det(g) \, g^{-1} ( \contra{v} \times \contra{w} )$,
		where $\contra{v}, \contra{w}$ denote the contravariant components.
		Comparing components yields the claim.
	\end{proof}
\end{proposition}

At the discrete level we use the \threeD complex \cref{eq: de rham 3D} to seek the velocity and mass flux through their flux representations in $\V^2(\mathcal{T})$, the depth in $\V^3(\mathcal{T})$ and the potential vorticity through its covariant representation $\covar{q}_h \in \V^1(\mathcal{T})$.
Writing $\covar{F}_h = g \flux{F}_h / \sqrt{g}$ for the covariant representation of the mass flux, 
and recalling \cref{lem: cross product,eq: mass pairings} means we replace the Coriolis term in \cref{eq: swe weak u} with
\begin{align}
	\int_{\thechart} \left( \covar{q}^{\mathrm{u}}_h \times \covar{F}_h \right) \cdot g \, \flux{\psi} \, \frac{1}{\sqrt{g}} \, \mathrm{d}\chart{x} ,
	\qquad
	\covar{q}^{\mathrm{u}}_h = \covar{q}_h - \vartheta \left( \partial_t \covar{q}_h - \frac{1}{\sqrt{g}} \, \flux{u}_h \times ( \nabla \times \covar{q}_h ) \right) ,
	\label{eq: swe 3d coriolis}
\end{align}
where the upwinding term is the covariant representation of $\surfvec{u} \times ( \surfcurl \surfvec{q} )$.
The potential vorticity diagnostic follows from the curl identity \cref{eq: ibp curl} and the boundary conditions:
find $\covar{q}_h \in \V^1(\mathcal{T})$ such that
\begin{align}
	\int_{\thechart} \varphi_h \, \covar{q}_h \cdot g^{-1} \covar{\omega} \, \sqrt{g} \, \mathrm{d}\chart{x}
	= \int_{\thechart} ( \nabla \times \covar{\omega} ) \cdot \frac{g \, \flux{u}_h}{\sqrt{g}} \, \mathrm{d}\chart{x}
	- \int_{\partial\thechart} ( \covar{n} \times \covar{\omega} ) \cdot \frac{g \, \flux{u}_h}{\sqrt{g}} \, \mathrm{d}s
	+ \int_{\thechart} \covar{\coriolis} \cdot g^{-1} \covar{\omega} \, \sqrt{g} \, \mathrm{d}\chart{x} ,
	\label{eq: swe 3d vorticity}
\end{align}
for all $\covar{\omega} \in \V^1(\mathcal{T})$, where $\covar{\coriolis}$ denotes the covariant representation of the radial Coriolis field with normal component $\surf{\coriolis}$.
All remaining terms of \cref{eq: swe weak,eq: swe diagnostics weak} carry over verbatim with the \threeD pullback indexing.
Mass conservation holds as in \cref{prop: mass conservation}. 
For energy conservation, taking $\flux{\psi} = \flux{F}_h$ in the momentum equation annihilates the Coriolis term pointwise, since $( \covar{q}^{\mathrm{u}}_h \times \covar{F}_h ) \cdot g \flux{F}_h / \sqrt{g} = ( \covar{q}^{\mathrm{u}}_h \times \covar{F}_h ) \cdot \covar{F}_h = 0$, and the argument of \cref{prop: energy conservation} applies.

\subsubsection{Time integration}
\label{sec: time integration}

The semi-discrete systems are integrated with the \ac{ssprk3} \cite{Shu1988}.
At every stage, the prognostic variables are advanced in time, and the diagnostic variables updated via  \cref{eq: swe diagnostics weak}. The time derivative of the potential vorticity in \cref{eq: supg,eq: swe 3d coriolis} is approximated at the final \RK stage.
The time step scales linearly with the mesh size, and the stabilisation parameter is $\vartheta = \Delta t / 2$. Specific parameter values are reported for each experiment in \cref{sec: results}.

%%%%%%%%%%%%%%%%%%%%%%%%%%%%%%%%%%%%%%%%%%%%%%%%%%%%%%%%%%%%%%%%%%%%%%%%%%%%%%
\subsection{Surface Stokes problem}
\label{sec: stokes}

In the typical  extrinsic setting with piecewise-polynomial geometry,
the surface Stokes problem is a recurring challenge 
since tangency and $H^1$-conformity cannot be enforced simultaneously \cite{Hansbo2020,Gross2018,Lederer2020,DemlowNeilan2026}.
In an intrinsic framework the geometry is exact, so tangency holds by construction. 
We therefore discretise the surface Stokes problem using the grad-conforming vector-valued spaces developed in \cref{sec: grad conforming} above, 
where $H^1$-conformity is lost only along chart interfaces as per \cref{rem: tangent H1}.

Given a viscosity $\nu > 0$ and a reaction coefficient $\alpha \geq 0$, the surface Stokes problem reads for $n=2$: find a tangent velocity $\surfvec{u}$ and a pressure $\surf{\varphi}$ such that
\begin{align}
	- \nu \left( \surfgrad ( \surfdiv \surfvec{u} ) 
	+ \surfgrad^{\perp} ( \surfgrad^{\perp} \cdot \, \surfvec{u} ) \right)
	+ \alpha \surfvec{u}
	+ \surfgrad \surf{\varphi} = \surfvec{f} ,
	\qquad
	\surfdiv \surfvec{u} = 0
	\qquad \text{on } \thesurface .
	\label{eq: surface stokes}
\end{align}
The contravariant components of velocity are represented
in the grad-conforming space 
$\mathbb{W}(\mathcal{T})$ 
developed in \cref{sec: grad conforming}, while the pressure is represented through its pullback in a scalar Lagrangian space $Q(\mathcal{T})$, which is type $\V^0(\mathcal{T})$ of order compatible with the velocity. 
Then,
the intrinsic \ac{fe} formulation of \cref{eq: surface stokes} is: find $\contra{u}_h \in \mathbb{W}(\mathcal{T})$, $\varphi_h \in Q(\mathcal{T})$ such that
\begin{subequations}
	\label{eq: surface stokes weak}
	\begin{align}
		\nu \int_{\thechart} \left( \divg \contra{u}_h \, \divg \contra{\psi}
		+ \divgperp \contra{u}_h \, \divgperp \contra{\psi} \right) \sqrt{g} \, \mathrm{d}\chart{x}
		&+ \alpha \int_{\thechart} \contra{u}_h \cdot g \, \contra{\psi} \, \sqrt{g} \, \mathrm{d}\chart{x}
		- \int_{\thechart} \varphi_h \, \divg \contra{\psi} \, \sqrt{g} \, \mathrm{d}\chart{x}
		\nonumber \\
		&= \int_{\thechart} \contra{f} \cdot g \, \contra{\psi} \, \sqrt{g} \, \mathrm{d}\chart{x} ,
		\\
		\int_{\thechart} \omega \, \divg \contra{u}_h \, \sqrt{g} \, \mathrm{d}\chart{x}
		&= 0 ,
	\end{align}
\end{subequations}
for all $\contra{\psi} \in \mathbb{W}(\mathcal{T})$, $\omega \in Q(\mathcal{T})$, where $\contra{f}$ collects the contravariant components of $\surfvec{f}$. 
Every term in \cref{eq: surface stokes weak} is expressed intrinsically, and the velocity is exactly tangent, with inter-chart continuity properties of \cref{prop: grad conforming}.

\begin{remark}
	\label{rem: stokes variants}
	The viscous term in \cref{eq: surface stokes} is the Hodge Laplacian of the \twoD complex.
	This differs by zeroth-order terms from the Bochner and deformation (Boussinesq--Scriven) variants used in the surface Stokes literature \cite{Gross2018,DemlowNeilan2026}.
	These variants are recovered on the sphere for divergence-free fields, by adjusting the reaction coefficient $\alpha$.	
\end{remark}

\section{Numerical experiments}
\label{sec: results}

In this section, we present numerical experiments for the model problems outlined above, and quantify the accuracy and computational benefits of the intrinsic framework.
All experiments are conducted on the cubed sphere manifold, whose atlas, panel partition and non-polynomial metric are provided in \cref{app: cubed sphere}.
The framework is implemented in the Julia package \GridapGeoscience \cite{GridapGeo,Tambyah2026_GridapGeo}, which builds on the \gridap ecosystem \cite{Badia2020,Verdugo2022}. 
Large-scale simulations use \GridapDistributed \cite{Badia2022}, \GridapPest \cite{GridapP4est}, \GridapPETSc \cite{GridapPETSc} and \GridapSolvers \cite{Manyer2024} on the Gadi supercomputer at NCI Australia.
The source code is available on Zenodo \cite{Zenodo_cubedsphere}.

In all simulations, $\mathcal{T}^{\ell}$ denotes the atlas triangulation with $N = 6 \times 2^{n\ell}$ cells, where $\ell = 0, 1, \dots, L$ is the refinement level, so that $\ell \rightarrow L$ corresponds to $h \rightarrow 0$.
Topological dimension $n = 2$ represents a spherical surface of radius $\radius$ without boundary, while $n = 3$ corresponds to an atmospheric shell of inner radius $\radius$ and thickness $\thickness$, with radial boundaries.
Ambient points $\surf{x} = (\surf{x}^1, \surf{x}^2, \surf{x}^3)$ are related to the longitude $\lambda = \arctan( \surf{x}^2 / \surf{x}^1 )$ and latitude $\theta = \arcsin( \surf{x}^3 / \surf{r} )$, where $\surf{r}^2 = (\surf{x}^1)^2 + (\surf{x}^2)^2 + (\surf{x}^3)^2$. 
Exact solutions and initial conditions specified on $\thesurface$ are converted into parametric representations using the pullbacks \cref{eq: pullbacks 3D,eq: pullbacks 2D}. Boundary data is prepared as described in \cref{rem: boundary data}.
%%%%%%%%%%%%%%%%%%%%%%%%%%%%%%%%%%%%%%%%%%%%%%%%%%%%%%%%%%%%%%%%%%%%%%%%%%%%%%
\subsection{Geometric accuracy under numerical quadrature}
\label{sec: quadrature study}

As per \cref{rem: exact mimetics}, the only geometry-related error of the intrinsic framework is the numerical integration of the metric-weighted terms.
So we now empirically
determine the degree of quadrature required to capture the non-polynomial cubed sphere metric.
We consider the primal form of the scalar Poisson problem \cref{eq: poisson mixed}: find $\varphi_h \in \V^0(\mathcal{T})$, subject to Dirichlet boundary conditions, such that
\begin{align}
	\int_{\thechart} \nabla \xi \cdot g^{-1} \nabla \varphi_h \, \sqrt{g} \, \mathrm{d}\chart{x}
	= \int_{\thechart} f \, \xi \, \sqrt{g} \, \mathrm{d}\chart{x} ,
	\qquad \forall \xi \in \V^0(\mathcal{T}) ,
	\label{eq: poisson primal}
\end{align}
which follows from the covariant mass pairing \cref{eq: mass pairings} and the integration by parts \cref{eq: ibp grad div}.
The solution is manufactured using the second-order polynomial pullback
\begin{align}
	\varphi_{\mathrm{ex}} = ( \pi/4 + x^1 ) ( \pi/4 - x^1 ) + ( \pi/4 + x^2 ) ( \pi/4 - x^2 ) ,
	\label{eq: poisson manufactured}
\end{align}
posed on a single panel of the cubed sphere in two and three dimensions, so that $\varphi_{\mathrm{ex}}$ belongs to the discrete space of Lagrange elements of order $2$. 
The forcing $f$ is a non-polynomial function that is manufactured via the intrinsic operators and the metric.
We use tensor-product Gauss--Legendre quadratures of degree $q$ on quadrilateral and hexahedral cells, where $q$, a positive odd integer, is the maximum polynomial degree integrated exactly.
On a flat domain, where $g$ is the identity, the exact solution is recovered to machine precision for $q \geq 5$ by the standard C\'ea and Galerkin orthogonality arguments, which assume exact integration.
On the cubed sphere, the departure from machine precision isolates the quadrature error to integrating the non-polynomial metric.

\begin{figure}[h!]
	\centering
	\includegraphics[width=\textwidth]{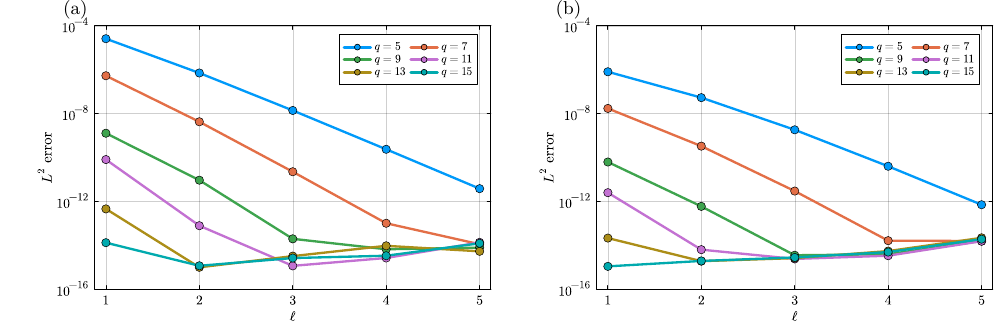}
	\caption{Convergence of the scalar Poisson problem in primal form \cref{eq: poisson primal} showing the $L^2$ error between exact and numerical solutions for (a) $n = 2$ and (b) $n = 3$, on a single panel with $2^{n\ell}$ cells.
	Colours indicate the quadrature degree $q$, for Lagrange elements of order $2$.
	Parameters: $\radius = 1$, $\thickness = 0.19$.}
	\label{fig: quadrature}
\end{figure}

\Cref{fig: quadrature} reports the $L^2$ norm of the error between exact and numerical solutions for Lagrange elements of order $2$, for varying quadrature degree $q$ and refinement level $\ell$.
As the norm of the exact solution is $\mathcal{O}(1)$, absolute and relative errors coincide up to a constant.
To ensure reliability, the reported errors are integrated with a conservative quadrature degree of $31$, and  the linear systems are solved with the UMFPACK sparse direct solver.
At the coarsest level, $q = 5$ captures the exact solution with $\mathcal{O}(10^{-4})$ error, which decreases to $\mathcal{O}(10^{-14})$ for $q = 15$. For fixed $q$, the error decreases towards machine precision as $\ell$ increases. For example, the $\mathcal{O}(10^{-14})$ accuracy of $q = 15$ at $\ell = 1$ is obtained with $q = 7$ at the finest level.
The integration accuracy can thus be increased either by raising $q$ at fixed $\ell$ or by refining the mesh at fixed $q$.
For all sufficiently large $q$ and $\ell$, the error settles at $\mathcal{O}(10^{-14})$. 
In this regime, 
the slight error growth with $\ell$ is consistent with the growth of the condition number of the \ac{fe} matrix.
These results confirm that the geometrical map is captured to machine precision via integration in the parametric space, meaning the intrinsic formulation is free of geometric consistency errors for a sufficient quadrature degree, as anticipated in \cref{sec: introduction}.

%%%%%%%%%%%%%%%%%%%%%%%%%%%%%%%%%%%%%%%%%%%%%%%%%%%%%%%%%%%%%%%%%%%%%%%%%%%%%%
\subsection{Convergence of the Hodge Laplacian problems}
\label{sec: hodge convergence}

\begin{figure}[h!]
	\centering
	\includegraphics[width=\textwidth]{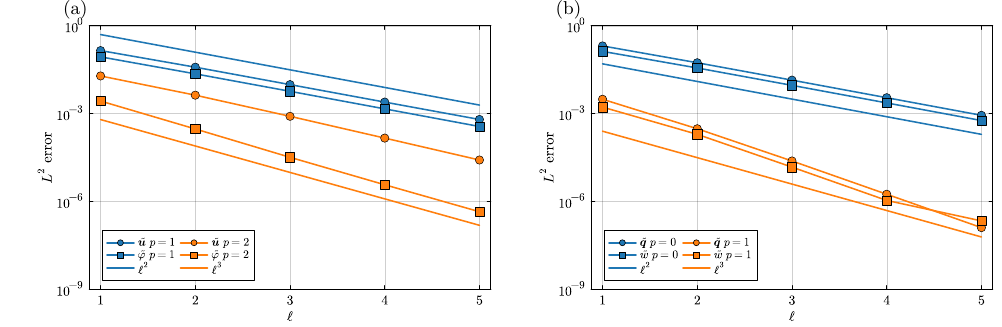}
	\caption{Convergence of the mixed Hodge Laplacian problems showing the $L^2$ error between exact and numerical solutions for (a) the scalar Poisson problem \cref{eq: poisson weak} and (b) the vector Poisson problem \cref{eq: vector poisson weak}, with $6 \times 2^{n\ell}$ cells, where $p$ is the index of the discrete complex \cref{eq: discrete complex}.
		Parameters: $n = 3$, $\radius = 1$, $\thickness = 0.19$.}
	\label{fig: hodge}
\end{figure}

The method of manufactured solutions is used to assess the convergence of the mixed formulations \cref{eq: poisson weak,eq: vector poisson weak} on an atmospheric shell for $n=3$.
The exact scalar solution is $\surf{\varphi}_{\mathrm{ex}} = \sin\theta$, and the exact vector solution $\surfvec{q}_{\mathrm{ex}}$ is the radial field of magnitude $2\sin\theta$. The corresponding auxiliary variables, forcing terms and boundary data are manufactured in parametric form, accordingly.
Convergence is measured in the $L^2$ norm, and the linear systems are solved with the MUMPS sparse direct solver \cite{GridapPETSc,MUMPS:1,MUMPS:2}.
The scalar problem \cref{eq: poisson weak} is discretised with $\flux{u}_h \in \V^2$ and $\varphi_h \in \V^3$, and the vector problem \cref{eq: vector poisson weak} with $\covar{q}_h \in \V^1$ and $w_h \in \V^0$. The value of $p$ reported in \cref{fig: hodge} below is the index of these spaces in the discrete complex \cref{eq: discrete complex}.
\Cref{fig: hodge} shows the expected $hp$-convergence for both the scalar and the vector Poisson problems.

%%%%%%%%%%%%%%%%%%%%%%%%%%%%%%%%%%%%%%%%%%%%%%%%%%%%%%%%%%%%%%%%%%%%%%%%%%%%%%
\subsection{Computational cost}
\label{sec: cost}

We now quantify the computational cost of the intrinsic pipeline discussed in \cref{sec: assembly} against the standard extrinsic one using the weak form of the Laplace--Beltrami operator. 
The pipelines differ in how the exact geometry enters the computation of the elemental stiffness integral $\int_{\thesurface} \surfgrad u \cdot \surfgrad v \, \mathrm{d}\thesurface$. 
Here \textit{extrinsic} refers to the standard assembly pipeline that relies on the Jacobian of the embedding. For the \twoD cubed sphere, this Jacobian is a $3 \times 2$ matrix, 
for which the pseudo-inverse is required to compute the surface gradient in the ambient space.
In contrast, the intrinsic framework introduced in \cref{sec: math form} above requires the metric to perform operations in the parametric space.
We consider a single quadrilateral cell $K$ of a chart in the  cubed sphere, which is obtained from the reference cell $\widehat{K}$ through the affine map $\phi : \widehat{K} \rightarrow K$. On $\widehat{K}$, we consider a tensor-product Gauss--Legendre quadrature of degree $q$ with points and weights $( \widehat{x}_i, \widehat{w}_i )$ as defined in \cref{sec: quadrature study} above \cite[Chapter 30]{ErnGuermond_2}.
In the extrinsic approach, the cell $\map \circ \phi ( \widehat{K} ) \subset \mathbb{R}^3$ is described using the Jacobian $J = \mathrm{d} ( \map \circ \phi ) \in \mathbb{R}^{3 \times 2}$, and the stiffness integral is approximated as 
\begin{align}
	\sum_{i} \left\{ \left( J^{+T} \, \widehat{\nabla} u \right) \cdot \left( J^{+T} \, \widehat{\nabla} v \right) \sqrt{ \det \left( J^T J \right) } \right\} ( \widehat{x}_i ) \; \widehat{w}_i ,
	\label{eq: extrinsic assembly}
\end{align}
where $(\cdot)^{+}$ denotes the Moore--Penrose pseudo-inverse, and $\widehat{\nabla}$ is the gradient in the reference cell.
The intrinsic pipeline evaluates the same integral entirely in the parametric space,
\begin{align}
	\sum_{i} \left\{ \widehat{\nabla} u \cdot \left( g^{-1} \, \widehat{\nabla} v \right) \sqrt{g} \right\} ( \widehat{x}_i ) \; w_i ,
	\label{eq: intrinsic assembly}
\end{align}
where $w_i$ absorbs the constant scaling of the affine map $\phi$, and $g^{-1} \sqrt{g}$ is evaluated using closed-form expressions \cref{eq: cubed sphere metric}. The assembly process in \cref{eq: intrinsic assembly} is that of the scalar Poisson problem \cref{eq: poisson primal} in \cref{sec: quadrature study}.

\Cref{tab: benchmark} reports the ratio of extrinsic \cref{eq: extrinsic assembly} to intrinsic \cref{eq: intrinsic assembly} cost indicators for $q = 11$, with number of floating point operations (column labelled `Flops'),
number of floating point operations per second (column labelled `Flops/sec.') and runtime (column labeled `Time') measured via GFlops.jl \cite{GFlops}.
As the \ac{fe} order, and hence the number of shape functions, increases, the extrinsic approach consistently requires nine times more operations than the intrinsic one, which correlates with longer runtimes.
The overhead is structural since \cref{eq: extrinsic assembly} maps \twoD reference quantities into the \threeD ambient space through $J$ and its pseudo-inverse at every quadrature point. 
In contrast, \cref{eq: intrinsic assembly} involves the $2 \times 2$ analytic metric such that every quantity is computed in the parametric space, as discussed in \cref{sec: assembly}.
The same conclusion extends to the other terms of the formulations of \cref{sec: models}. In particular, the rotated pairings in \cref{lem: rotation} that replace \threeD cross products by inexpensive \twoD operations.

\begin{table}[h!]
	\centering
	\begin{tabular}{cccc}
		\toprule
		\ac{fe} order & Flops/sec. & Flops & Time \\
		\midrule
		1 & 3.028 & 9 & 2.972 \\
		\midrule
		2 & 2.433 & 9 & 3.699 \\
		\midrule
		3 & 2.248 & 9 & 4.003 \\
		\bottomrule
	\end{tabular}
	\vspace{0.5em}
	\caption{Ratio of extrinsic \cref{eq: extrinsic assembly} to intrinsic \cref{eq: intrinsic assembly} cost indicators for the weak form of the Laplace--Beltrami operator on a quadrilateral cell, with quadrature degree $q = 11$.}
	\label{tab: benchmark}
\end{table}
%%%%%%%%%%%%%%%%%%%%%%%%%%%%%%%%%%%%%%%%%%%%%%%%%%%%%%%%%%%%%%%%%%%%%%%%%%%%%%
\subsection{Rotating shallow water equations}
\label{sec: swe experiments}

The two- and \threeD semi-discrete systems of \cref{sec: swe} are integrated in time as described in \cref{sec: time integration}, 
where the mass systems \cref{eq: swe diagnostics weak} are solved using the conjugate gradient method \cite{Manyer2024}.
For $n = 2$ the prognostic variables are $\flux{u}_h \in \V^1$ and $\varphi_h \in \V^2$, while for $n = 3$ they are $\flux{u}_h \in \V^2$ and $\varphi_h \in \V^3$, with the potential vorticity $\covar{q}_h \in \V^1$. As above, the reported $p$ is the index of these spaces in the discrete complex \cref{eq: discrete complex}.
The time step is $\Delta t = C \, \Delta x$ with $C = \mathrm{CFL} / (p+1)^2 \sqrt{\gravity H_0}$, where $\Delta x = \sqrt{4 \pi \radius^2 / N}$ is the nominal horizontal resolution, $H_0$ is the mean depth, and $\vartheta = \Delta t / 2$.

\subsubsection{Williamson 2: steady zonal flow}

\begin{figure}[h!]
	\centering
	\includegraphics[width=\textwidth]{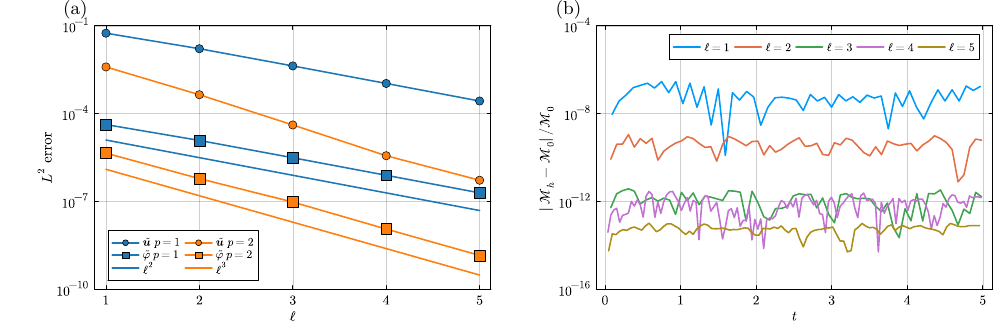}
	\caption{Williamson 2 test case: (a) the relative $L^2$ error between initial and final solutions, and (b) the mass conservation error for $p = 1$, with time in days, where $p$ is the index of the discrete complex \cref{eq: discrete complex}.
	Parameters: $u_0 = 0.083$, $\gravity = 289.49$, $H_0 = 4.7 \times 10^{-4}$, $\surf{\topography} = 0$, $\mathrm{CFL} = 0.1$, $\radius = 1$, $6 \times 2^{2\ell}$ cells.}
	\label{fig: W2}
\end{figure}

Convergence of the \twoD system \cref{eq: swe weak,eq: swe diagnostics weak}  is assessed using the Williamson 2 test case, which is a steady state solution of the shallow water equations \cite{Williamson1992}.
Simulations are conducted in a non-dimensional frame corresponding to the length scale $\widehat{L} = 6371220~\mathrm{m}$ and time scale $\widehat{T} = ( 7.292 \times 10^{-5} )^{-1}~\mathrm{s}$.
The initial conditions are
\begin{align}
	\surf{\varphi} = H_0 - \frac{u_0}{\gravity} \left( 1 + \frac{1}{2} u_0 \right) \sin^2\theta - \surf{\topography} ,
	\qquad
	\surfvec{u} = u_0 \cos\theta \; \surfvec{e}_{\lambda} ,
	\label{eq: W2 IC}
\end{align}
where $\surfvec{e}_{\lambda}$ is the unit vector in the zonal direction, with parameters $u_0 = 0.083$, $H_0 = 4.7 \times 10^{-4}$, $\gravity = 289.49$, $\surf{\topography} = 0$ and Coriolis parameter $\surf{\coriolis} = 2 \sin\theta$ \cite{Williamson1992}.
Similar to other studies \cite{Rognes2013}, convergence is assessed using the relative $L^2$ error between the initial condition and the solution after five days. 
\Cref{fig: W2}(a) shows the expected second- and third-order convergence for $p = 1$ and $p = 2$ elements,
where spatial and temporal resolution reduce simultaneously at a fixed $\mathrm{CFL}$.
\Cref{fig: W2}(b) reports the normalised mass conservation error over time for $p = 1$ elements, which decreases with the refinement level and reaches round-off levels at the finest resolutions.
Semi-discrete mass conservation is exact by \cref{prop: mass conservation}. 
The residual error observed at coarse resolutions is not a property of the formulation but of its algebraic solution
since the relative tolerance of the solver is $10^{-8}$.

\subsubsection{Williamson 5: zonal flow over an isolated mountain}

The Williamson 5 test case describes zonal flow over an isolated mountain \cite{Williamson1992}. This test case is used to confirm that the \threeD system of \cref{sec: swe 3d} reproduces the \twoD dynamics at lowest order.
The initial condition is \cref{eq: W2 IC} with $u_0 = 0.043$, $H_0 = 9.4 \times 10^{-4}$, $\gravity = 289.49$ and the topography
\begin{align}
	\surf{\topography} = b_0 \left( 1 - \nu \eta^{-1} \right) ,
	\qquad
	\nu = \min \left[ \eta , \, \sqrt{( \lambda - \lambda_c )^2 + ( \theta - \theta_c )^2} \right] ,
	\label{eq: W5 topography}
\end{align}
where $b_0 = 3 \times 10^{-4}$, $\eta = \pi/9$ and the mountain centred at $( \lambda_c, \theta_c ) = ( 0, \pi/6 )$.
The duration of $20$ days is sufficient for a Rossby wave to develop and meander around the mountain.

\begin{figure}[h!]
	\centering
	\includegraphics[width=0.9\textwidth]{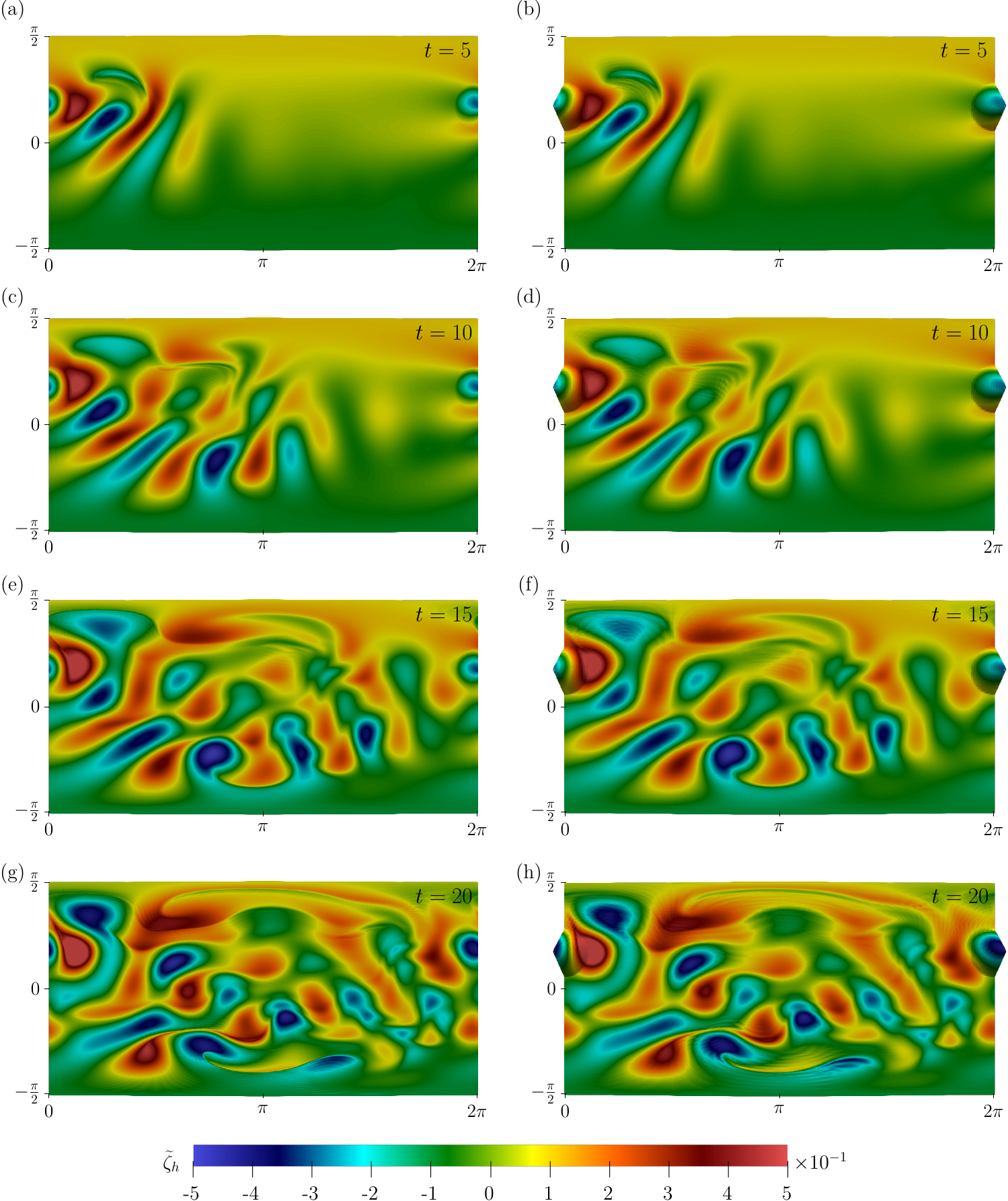}
	\caption{Williamson 5 test case: (a,c,e,g) \twoD relative vorticity, and (b,d,f,h) radial component of the \threeD relative vorticity (aspect ratio exaggerated for visualisation), at $t = 5, 10, 15, 20$ days.
		Parameters: $b_0 = 3 \times 10^{-4}$, $u_0 = 0.043$, $H_0 = 9.4 \times 10^{-4}$, $\gravity = 289.49$, $\eta = \pi/9$, $( \lambda_c, \theta_c ) = ( 0, \pi/6 )$, $\ell = 6$, $\mathrm{CFL} = 0.1$, $\radius = 1$, $\thickness = H_0$.}
	\label{fig: W5}
\end{figure}

\begin{figure}[h!]
	\centering
	\includegraphics[width=\textwidth]{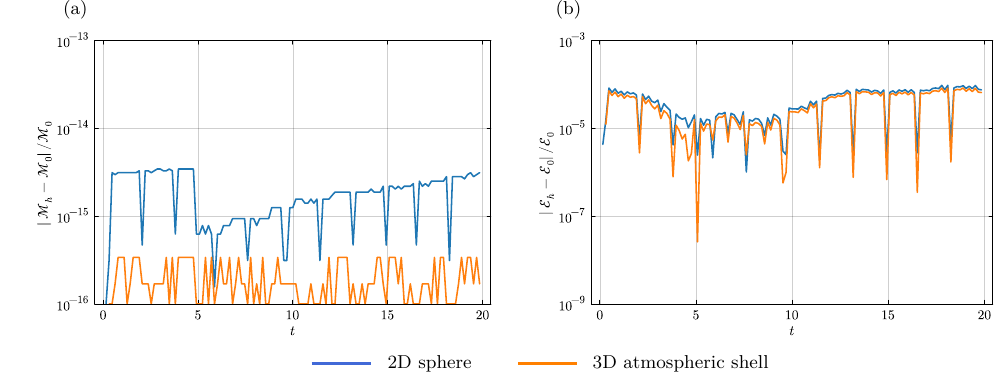}
	\caption{Williamson 5 test case showing normalised conservation errors of (a) mass and (b) energy, comparing the \twoD ($p = 1$) and \threeD ($p = 0$) systems, with time in days and parameters as in \cref{fig: W5}.}
	\label{fig: W5 conservation}
\end{figure}

In two dimensions, we solve \cref{eq: swe weak,eq: swe diagnostics weak} with $p = 1$ elements at refinement level $\ell = 6$, with $\surf{\coriolis} = 2 \sin\theta$.
\Cref{fig: W5}(a,c,e,g) shows the evolution of the relative vorticity $\surf{\zeta}_h = \surfgrad^{\perp} \cdot \, \surfvec{u}_h$,  and the expected Rossby wave meander.
The potential vorticity upwinding \cref{eq: supg} suppresses grid-scale oscillations, and turbulent dynamics are resolved over long time scales.
\Cref{fig: W5 conservation}(a) shows that discrete mass is conserved to machine precision, while the discrete energy \cref{eq: energy} is conserved to $\mathcal{O}(10^{-4})$ over the simulation. 
The semi-discrete conservation of \cref{prop: energy conservation} is perturbed only by the explicit \ac{ssprk3} time discretisation.
In situations where fully discrete energy conservation is required, an implicit, conservative time integrator can be employed \cite{Cohen2011,BauerCotter2018,Wimmer2020}.

In three dimensions,  the Coriolis force has normal component $2 \sin\theta$, and we solve the system of \cref{sec: swe 3d} 
with $p=0$ elements and the same time step as the \twoD case.
The computational domain consists of two atmospheric shells at horizontal refinement level $\ell = 6$, which are connected via a single cell in the radial direction. 
The parametric coordinates of the bottom and top shells are $(x^1, x^2, 0)$ and $(x^1, x^2, 1)$, respectively.
The effect of topography is included as a \emph{perturbation of the geometry} such that the bottom shell is displaced radially by $\surf{\topography}$. That is, the bottom shell sits at radius $\radius + \surf{\topography}$, while the top shell is flat at radius $\radius + \thickness$ where $\thickness = H_0$ reflects the depth of the atmosphere.
\Cref{fig: W5}(b,d,f,h) shows the radial component of the relative vorticity in three dimensions, which is comparable to the \twoD solution. Since $\thickness \ll \radius$, the aspect ratio  in \Cref{fig: W5}(b,d,f,h) is exaggerated for visualisation.
The minor discrepancies between the \twoD and \threeD conservation errors in \cref{fig: W5 conservation} is likely due to the different \ac{fe} orders employed.
This experiment demonstrates that the intrinsic framework is agnostic to the choice of geometrical map, and extends naturally to atmospheric configurations with orography.

%%%%%%%%%%%%%%%%%%%%%%%%%%%%%%%%%%%%%%%%%%%%%%%%%%%%%%%%%%%%%%%%%%%%%%%%%%%%%%
 
\subsection{Surface Stokes problem}
\label{sec: stokes experiments}
 
\newcommand{\vertiii}[1]{{\left\vert\kern-0.05ex\left\vert\kern-0.05ex\left\vert #1 \right\vert\kern-0.05ex\right\vert\kern-0.05ex\right\vert}}

We now assess the nodal construction developed in \cref{sec: nodal construction} above, whose vector-valued fields are exactly tangent and inter-chart continuous at the nodes, through the convergence of the surface Stokes problem \cref{eq: surface stokes weak}.
The exact solutions are  
$\surfvec{u}_{\mathrm{ex}}=-\sin \lambda \surfvec{e}_{\lambda}$ 
and $\surf{\varphi}_{\mathrm{ex}} = \sin \theta $, for which  corresponding forcing terms are manufactured accordingly. 
We use the Taylor--Hood \ac{fe} pair \cite{DemlowNeilan2026}, so that the velocity space $\mathbb{W}(\mathcal{T})$ is of order $p+1$ and the pressure space $Q(\mathcal{T})$ of order $p$.  The linear system is solved via the  MUMPS sparse direct solver \cite{GridapPETSc,MUMPS:1,MUMPS:2}, and all numerical solutions and reported errors are integrated with sufficiently conservative quadrature degrees.

\begin{figure}[h!]
	\centering
	\includegraphics[width=\textwidth]{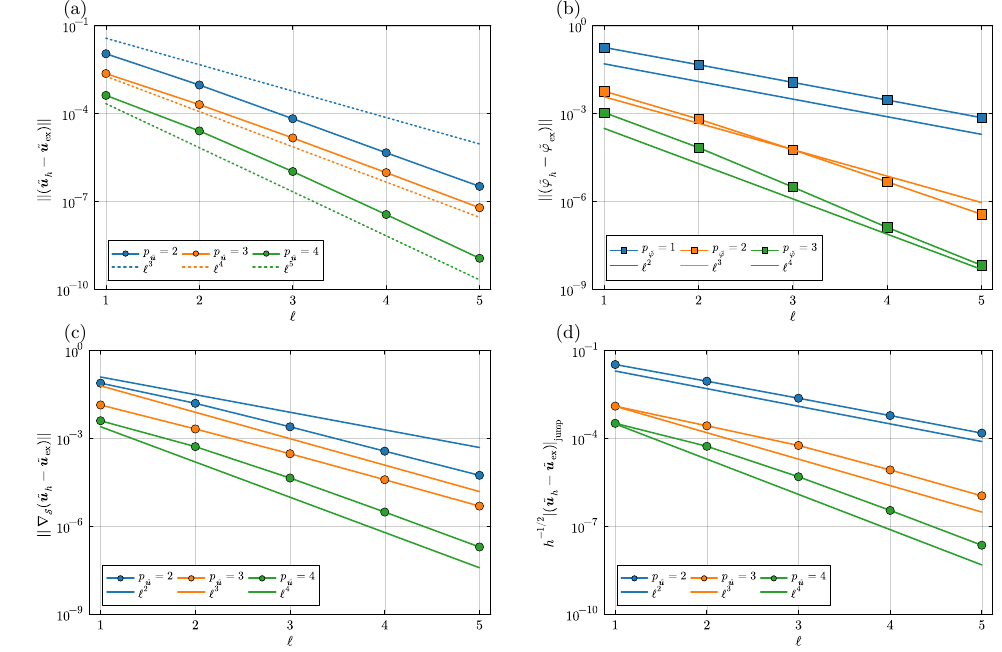}
	\caption{Convergence of surface Stokes problem showing the error between exact and numerical solutions for the (a,b) $L^2$ norm for velocity and pressure,  (c) gradient norm for velocity, and (d) jump norm for velocity \cref{eq: jumps}, with $6 \times 2^{n\ell}$ cells, where $p_{\surfvec{u}} = p+1$ and $p_{\surf{\varphi}} = p$ denote the orders of the velocity and pressure spaces. 
		Parameters: $n = 2$, $\radius = 1$,  $\nu=1$, $\alpha=1$.}
	\label{fig: stokes v2}
\end{figure}

Convergence is assessed using the following discontinuous Galerkin-like norm \cite{Cockburn2009} that accounts for the jump between chart interfaces: 
\begin{align}
	\vertiii{\surfvec{u}}_h^2 = 
	\alpha || \surfvec{u}||_{L^2(\thesurface)}^2
	+
	\sum_{\surf{K}}|| \surfgrad \surfvec{u}||_{L^2(\surf{K})}^2
	+
	\sum_{\surf{e} } h_{\surf{e}}^{-1} || \llbracket\surfvec{u} \otimes \surfvec{n}\rrbracket  ||_{L^2(\surf{e})}^2,
	\label{eq: jumps}
\end{align}
where $	\llbracket \surfvec{u} \otimes \surfvec{n} \rrbracket= \surfvec{u}^+\otimes \surfvec{n}^+ 
+ \surfvec{u}^- \otimes \surfvec{n}^-$
is the jump across interior edges $\surf{e}$, and $\surfvec{n}$ is the normal vector to the edge in the tangent space of the manifold.
In \cref{eq: jumps}, the first $L^2$ norm is expected to converge at rate $h^{p+2}$,	while the second gradient norm and the third jump norm are both expected to converge at rate $h^{p+1}$.
\Cref{fig: stokes v2} reports the convergence of each norm in \cref{eq: jumps} individually.
\Cref{fig: stokes v2}(a,b) shows the expected $hp$-convergence of the $L^2$ norm of the velocity and scalar solutions, respectively, with the velocity converging at rate $h^{p+2}$. 
For the velocity, \cref{fig: stokes v2}(c,d) illustrates that the gradient norm and the jump norm scaled by $h^{-1/2}$, respectively, converge at rate $h^{p+1}$. 
The jump contribution therefore converges at the same rate as the broken gradient, so the interface non-conformity does not degrade the rate measured in \cref{eq: jumps}. A bound on the interface jump is given in \cref{rem: accuracy}.
A detailed stability and convergence analysis for the surface Stokes and $H^1$ vector Laplacian problems is left for future work.

\section{Conclusions and future work}
\label{sec: conclusion}

We develop an intrinsic \ac{fe} framework for scalar- and vector-valued \acp{pde} on manifolds described by general atlases of charts.
The geometry enters the formulation only through the metric tensor,
which is evaluated at quadrature points.
The mimetic structure of the discrete de Rham complex holds  independently of the accuracy of numerical quadrature, 
since 
differential operators, duality pairings between complementary representations, annihilation properties, and integration by parts identities  
are metric-free.
\ac{fe} spaces of arbitrary order are obtained from standard Euclidean elements in the parametric domains.
For the spaces of the de Rham complex, the global spaces are conforming and inter-chart gluing is purely topological, since the canonical degrees of freedom are chart-invariant.
For grad-conforming vector-valued spaces,
we propose a nodal change of basis construction that yields exactly tangent discrete vector fields, where inter-chart continuity holds at the nodes, with interface jumps of the order of the approximation error.
Numerical experiments on the cubed sphere manifold confirm the  non-polynomial geometry is captured to machine precision for a sufficient quadrature degree, 
and that intrinsic assembly requires an order of magnitude fewer operations than its extrinsic counterpart.  
Case studies show the Hodge Laplacian and surface Stokes problems converge at the expected rates, 
and that shallow water discretisations, with orography included as a perturbation of the geometry, conserve mass exactly and energy up to the time discretisation error, 
and sustain nonlinear dynamics over long time scales.

There are several possible extensions of this work.
The framework is agnostic to the geometrical map and to the atlas, so application to non-spherical and non-zero genus manifolds is feasible provided a discrete representation of the corresponding parametric space. 
Similarly, extension to terrain-following \threeD geometries provide a natural path towards the compressible Euler equations of operational atmospheric models \cite{Wood2013,Melvin2024_orography}.
Further extensions include an a priori error analysis accounting for the quadrature of the metric terms that relates to geometric variational crime theory \cite{HolstStern2012,Holst2023}, 
and generalisation of the vector calculus--based presentation via \ac{fe} exterior calculus \cite{Arnold2018}.
As discussed in \cref{sec: stokes}, 
our novel construction of grad-conforming vector-valued spaces is amenable to surface Stokes discretisations with exactly tangential velocities {where inter-chart continuity holds at the nodes}.
The corresponding analysis, including a rigorous convergence theory for the nodal construction, is left for future consideration.
Finally, Stokes complexes, whose leading spaces require $C^1$ elements, are beyond the scope of this work \cite{Falk2013,Chen2024}.

\section*{Acknowledgements}
This research is supported by the Commonwealth of Australia as represented by the Defence Science and Technology Group of the Department of Defence.
This research is also funded by the Australian Government through the Australian Research Council (project numbers DP210103092 and DP220103160).
This research was undertaken with the assistance of resources from the Monash University and National Computational Infrastructure (NCI Australia) allocation scheme. NCI is an NCRIS-enabled capability supported by the Australian Government.
This work is also supported by computational resources provided by the Australian Government through NCI under the National Computational Merit Allocation Scheme (NCMAS) and ANU Merit Allocation Scheme.

\subsection*{Author contribution}
All authors contributed equally to the design of the study.
Tamara A. Tambyah and Santiago Badia drafted the article, and Tamara A. Tambyah performed the numerical simulations with the assistance of Alberto F. Mart\'in.
All authors gave approval for publication.

\subsection*{Declaration of competing interest}
The authors declare they have no competing interests.

\subsection*{Data availability}
This article does not contain any additional data.
The source code for this study is available on Zenodo \cite{Zenodo_cubedsphere}.

\appendix

\section{The cubed sphere atlas}
\label{app: cubed sphere}

The cubed sphere manifold \cite{Ronchi1996} is obtained by projecting the six faces of a cube onto the sphere of radius $\radius$.
The six charts of the atlas, referred to as \emph{panels}, consist of the enlarged\footnote{The enlargement makes neighbouring panels genuinely overlap, so that the transition maps are defined on neighbourhoods of the panel interfaces, as required in \cref{sec: manifolds}. The implementation \cite{Tambyah2026_GridapGeo} stores the partitions $P$, which only touch along their closed edges, rather than the enlarged chart domains.} open squares $\thechart = ( -\pi/4 - \epsilon, \pi/4 + \epsilon )^2$, with $0 < \epsilon < \pi/4$, 
while the partition domains are the closed squares $P = [ -\pi/4, \pi/4 ]^2$.
The metric tensor and measure is
\begin{align}
	g = \frac{ \radius^2 }{ \rho^4 \cos^2 x^1 \cos^2 x^2 }
	\begin{pmatrix}
		1 + \tan^2 x^1 & - \tan x^1 \tan x^2 \\
		- \tan x^1 \tan x^2 & 1 + \tan^2 x^2
	\end{pmatrix} ,
	\qquad
	\sqrt{g} = \frac{ \radius^2 }{ \rho^3 \cos^2 x^1 \cos^2 x^2 } ,
	\label{eq: cubed sphere metric}
\end{align}
where $\rho = \sqrt{ 1 + \tan^2 x^1 + \tan^2 x^2 }$ and $x^{1,2}$ are parametric coordinates.
This geometry corresponds to the well-known equiangular projection, 
for which \cref{eq: cubed sphere metric} is identical on all six panels \cite{Nair2005}.
Other geometrical maps are also possible \cite{Rancic1996,Putman2007,Mcgregor2005,Giraldo2003}.

Extruding the \twoD cubed sphere manifold in the radial direction
by a constant thickness $\thickness \ll \radius$ represents an atmospheric shell, with metric and measure 
\begin{align}
	g_{3\mathrm{D}} =
	\begin{pmatrix}
		g \big|_{ \radius \rightarrow \radius + \thickness x^3 } & 0 \\
		0 & \thickness^2
	\end{pmatrix} ,
	\qquad
	\sqrt{ g_{3\mathrm{D}} } = \thickness \, \sqrt{g} \, \Big|_{ \radius \rightarrow \radius + \thickness x^3 } ,
	\label{eq: cubed sphere shell metric}
\end{align}
where $g|_{\radius \rightarrow \radius + \thickness x^3}$ denotes the \twoD metric \cref{eq: cubed sphere metric} evaluated with the perturbed radius  $\radius + \thickness x^3$. 
The orography configuration in \cref{sec: swe experiments} is the variable-thickness instance of an atmospheric shell with radial coordinate $\radius + \topography( x^1, x^2 ) + x^3 ( \thickness - \topography( x^1, x^2 ) )$. 
That is, the topography displaces the bottom shell while the top shell is flat.
In this case, the metric contains off-diagonal terms proportional to the topography gradient.

\printbibliography

\end{document}